\documentclass[reqno,a4paper]{amsart}

\usepackage{graphicx, geometry}
\usepackage{color} 
\usepackage[usenames,dvipsnames]{xcolor}
\usepackage{transparent}
\usepackage[utf8x]{inputenc}  
\usepackage{dsfont}
\usepackage[active]{srcltx}\usepackage[colorlinks,pdfpagelabels,pdfstartview=FitH,bookmarksopen=true,bookmarksnumbered=true,linkcolor=blue,plainpages=false,hypertexnames=false,citecolor=red]{hyperref}
\usepackage{cleveref}
\usepackage{mathrsfs}
\usepackage{verbatim, amssymb,amsmath, amsfonts, amsthm}
\usepackage{latexsym}
\usepackage[dvips]{epsfig}
\usepackage{float}
\usepackage[normalem]{ulem}

\allowdisplaybreaks 

\newtheorem{teorema}{Theorem}[section]
\newtheorem{proposizione}[teorema]{Proposition}
\newtheorem{lemma}[teorema]{Lemma}

\newtheorem{osservazione}[teorema]{Remark}
\newtheorem{definizione}[teorema]{Definition}

\numberwithin{equation}{section}


\newcommand{\R}{\mathbb{R}}
\renewcommand{\to}{\rightarrow}

\newcommand{\upm}{u_p^{(m)}}
\newcommand{\upjm}{u_{p_j}^{(m)}}

\newcommand{\upLm}{u_{p,L}^{(m)}}
\newcommand{\hpm}{h_{p,m}}

\newcommand{\hpLm}{h_{p,m,L}}

\newcommand{\vpjdmu}{v_{p_j}^{(2m-1)}}

\newcommand{\vpLdmu}{v_{p,L}^{(2m-1)}}
\newcommand{\vpn}{v_p^{(n)}}

\newcommand{\vp}{v_p}
\newcommand{\vpin}{v_{p,n,i}}

\newcommand{\vpntildei}{\tilde v_{p,n,i}}

\newcommand{\barvpin}{\bar{v}_{p,n,i}}

\newcommand{\wpn}{w_{p,n}}

\newcommand{\wpntildei}{\tilde w_{p,n,i}}

\newcommand{\mupn}{\mu_{p,n}}

\newcommand{\hpjm}{{h}_{p_j,m}}
\def\sideremark#1{\ifvmode\leavevmode\fi\vadjust{\vbox to0pt{\vss
			\hbox to 0pt{\hskip\hsize\hskip1em
				\vbox{\hsize1.8cm\tiny\raggedright\pretolerance10000
					\noindent #1\hfill}\hss}\vbox to15pt{\vfil}\vss}}}%

\newcommand{\virg}[1]{\textquotedblleft#1\textquotedblright}

\begin{document}
	\date{\today}
	\title[]{One-dimensional nodal solutions of the Lane--Emden problem in cylinders:\\ asymptotic behaviour and energy stability}
	
	\author[]{Francesca De Marchis, Lisa Mazzuoli}
	
	\address{Francesca De Marchis, University of Roma {\em Sapienza}, P.le Aldo Moro 5, 00185 Roma, Italy -- demarchis@mat.uniroma1.it}
	\address{Lisa Mazzuoli, University of Roma {\em Sapienza}, P.le Aldo Moro 5, 00185 Roma, Italy -- lisa.mazzuoli@uniroma1.it}
	
	
	
	\thanks{Research partially supported by PRIN Project 2022AKNSE4 \emph{Variational and Analytical aspects of Geometric PDE} and by Gruppo Nazionale per l’Analisi Matematica, la Probabilità e le loro
		Applicazioni (GNAMPA) of the Istituto Nazionale di Alta Matematica (INdAM). In particular the first author is supported by INdAM-GNAMPA Project $\#$E53C25002010001$\#$.}

	\begin{abstract}
		We consider the Lane--Emden problem $-\Delta u=|u|^{p-1}u$, $p>1$, in
		bounded domains of a half-cylinder $\omega\times(0,+\infty)$, with
		Dirichlet conditions on the relative boundary and Neumann conditions on
		the boundary of the cylinder. In flat cylinders $\omega\times(0,L)$ this
		problem admits, for every $p>1$, one-dimensional solutions with any
		prescribed number $m$ of nodal regions. We study the stability of these
		domains, as critical points of a related energy functional, with respect to volume-preserving deformations of the domain. We prove that, for $p$
		large, the flat cylinders are stable as soon as their height exceeds a
		threshold growing linearly in $p$, while, for $p$ close to $1$, both
		stability and instability occur, with explicit thresholds depending on
		$m$ and on the first nontrivial Neumann eigenvalue of $\omega$. These
		results rely on a detailed analysis of the one-dimensional profiles: we
		describe their asymptotic behaviour, and that of the first eigenvalue of
		the associated linearized operator, as $p\to1^+$ and $p\to+\infty$, we
		obtain expansions of the amplitude with explicit constants, and we prove
		global comparison results with respect to $p$, some of which are new
		also for positive solutions.
	\end{abstract}
	\maketitle
	
	\section{Introduction}
	In this paper we study one-dimensional sign-changing solutions of a relative
	Lane--Emden problem in cylindrical domains. By a one-dimensional
	solution we mean a function depending only on one Cartesian variable.
	More precisely, if $x=(x',x_N)\in\mathbb R^{N-1}\times\mathbb R$, we
	consider solutions of the form
	\[
	u(x)=u(x',x_N)=u(x_N).
	\]
	The natural domains for such functions are cylinders, in which the
	$x_N$-variable represents a distinguished direction. One-dimensional solutions arise in various PDE problems; for instance,
	they are related to the celebrated De Giorgi conjecture (see \cite{GhoussoubGui} and the references therein).
	One-dimensional solutions in cylinders
	can be viewed as natural counterparts of radial solutions in balls, although they have received less attention.
	
	In particular in this paper we consider the unbounded half-cylinder 
	\[
	\mathcal C_\omega:=\omega\times(0,+\infty)
	\]
	spanned by a smooth bounded
	domain $\omega\subset\mathbb R^{N-1}$, $N\geq2$. We denote by
	\[
	\Gamma_\Omega:=\partial\Omega\cap\mathcal C_\omega\qquad\text{and}\qquad \Gamma_{1,\Omega}:=\partial\Omega\setminus\overline\Gamma_\Omega,
	\]
	where we refer to $\Gamma_\Omega$, which is assumed to be smooth, as the relative, or free, boundary and we denote by $\nu$ the outer unit normal vector to $\Gamma_\Omega$.\\
	We focus on the Lane--Emden problem
	\begin{equation}\label{problem_Omega}
		\begin{cases}
			-\Delta u=|u|^{p-1}u & \text{in }\Omega,\\
			u=0 & \text{on }\Gamma_\Omega,\\
			\dfrac{\partial u}{\partial\nu}=0 & \text{on }\Gamma_{1,\Omega},
		\end{cases}
	\end{equation}
	where $p>1$ and $\nu$ is the outer unit normal vector on $\Gamma_{1,\Omega}$.
	
	We restrict our attention to domains $\Omega\subset\mathcal C_\omega$
	whose relative boundary is a Cartesian graph over $\omega$. More
	precisely, for $\varphi\in C^2(\bar\omega)$ we denote by $\Omega_\varphi$
	the hypograph in $\mathcal C_\omega$ of the positive function $e^{\varphi}$,
	namely
	\[
	\Omega_\varphi:=\bigl\{(x',x_N)\in\mathcal C_\omega:\ x_N<e^{\varphi(x')}\bigr\},
	\]
	and we write
	\[
	\Gamma_\varphi:=\Gamma_{\Omega_\varphi}=\bigl\{(x',x_N)\in\mathcal C_\omega:\ x_N=e^{\varphi(x')}\bigr\},
	\qquad
	\Gamma_{1,\varphi}:=\Gamma_{1,\Omega_\varphi}.
	\]
	From now on, problem \eqref{problem_Omega} is set in domains of
	this type.
	
	If $1<p<p_S$, where $p_S:=\frac{N+2}{N-2}$ if $N\geq3$ and
	$p_S:=+\infty$ if $N=2$, problem \eqref{problem_Omega} in $\Omega_\varphi$
	can be studied variationally in the space
	$H^1_0(\Omega_\varphi\cup\Gamma_{1,\varphi})$ of functions in
	$H^1(\Omega_\varphi)$ whose trace vanishes on $\Gamma_\varphi$: its weak
	solutions are the critical points of the energy functional
	\[
	J(u):=\frac12\int_{\Omega_\varphi}|\nabla u|^2\,dx-\frac1{p+1}\int_{\Omega_\varphi}|u|^{p+1}\,dx,
	\qquad u\in H^1_0(\Omega_\varphi\cup\Gamma_{1,\varphi}).
	\]
	A positive solution can be obtained by minimizing $J$ on the Nehari
	manifold
	\[
	\mathcal N:=\bigl\{u\in H^1_0(\Omega_\varphi\cup\Gamma_{1,\varphi})\setminus\{0\}:\ J'(u)[u]=0\bigr\}
	\]
	(see \cite{AmbrosettiRabinowitz, CastroCossioNeuberger}), while, arguing
	as in \cite{CastroCossioNeuberger, BartschWeth}, a least energy
	sign-changing solution can be obtained by minimizing $J$ on the nodal
	Nehari set
	\[
	\mathcal M:=\bigl\{u\in H^1_0(\Omega_\varphi\cup\Gamma_{1,\varphi}):\ u^+\neq0,\ u^-\neq0,\ J'(u)[u^+]=J'(u)[u^-]=0\bigr\},
	\]
	where $u^+:=\max\{u,0\}$ and $u^-:=\min\{u,0\}$, and it has exactly two
	nodal regions. Moreover, since the nonlinearity is odd, problem
	\eqref{problem_Omega} admits infinitely many solutions
	\cite{AmbrosettiRabinowitz}, but in general without a precise control
	of their nodal structure.
	
	In general, given a solution $u_{\Omega_\varphi}$ of \eqref{problem_Omega}
	in $\Omega_\varphi$, we may consider small deformations of $\Omega_\varphi$
	within the class of domains introduced above, namely
	$\Omega_{\varphi+t\psi}$ with $\psi\in C^2(\bar\omega)$ and $|t|$ small.
	The nondegeneracy of $u_{\Omega_\varphi}$ induces local uniqueness for
	solutions of \eqref{problem_Omega} in the deformed domains (see
	Section~\ref{section:stability} and \cite[Proposition 2.1]{AfonsoPacella}).
	Thus, for a nondegenerate solution $u_{\Omega_\varphi}$ of
	\eqref{problem_Omega}, an interesting question is to understand how the
	energy behaves with respect to variations of $\Omega_\varphi$ which
	preserve its volume.
	
	If $\Omega$ is the flat cylinder
	\[
	\Omega_{\varphi_L}=\omega\times(0,L),\qquad \varphi_L\equiv\log L,\quad L>0,
	\]
	then, for every $p>1$ and every integer $m\geq1$, a particular solution
	of \eqref{problem_Omega} is obtained by extending to $\Omega_{\varphi_L}$
	a solution of
	\begin{equation}\label{ODE_in_(0,L)}
		\begin{cases}
			-u''=|u|^{p-1}u & \text{in }(0,L),\\
			u(L)=0,\\
			u'(0)=0,
		\end{cases}
	\end{equation}
	having $m$ nodal regions. Such a solution is unique up to sign (see
	Proposition~\ref{prop:existence_uniqueness} and the even reflection
	described below), and we select the one satisfying
	\[
	(u_{p,L}^{(m)})'(L)<0.
	\]
	With a slight abuse of notation, we denote by the same symbol its
	one-dimensional extension to the cylinder,
	\begin{equation}\label{write_sol}
		u_{p,L}^{(m)}(x',x_N):=u_{p,L}^{(m)}(x_N),\qquad (x',x_N)\in\Omega_{\varphi_L},
	\end{equation}
	which is a solution of \eqref{problem_Omega} in $\Omega_{\varphi_L}$ with
	$m$ nodal regions.
	
	The
	nondegeneracy of $u_{p,L}^{(m)}$, instead, depends on $\omega$ and $L$:
	throughout the paper we work under sufficient conditions on $\omega$ and
	$L$ which guarantee it (see \eqref{nondeg_cond} below). In this case,
	for $\varphi$ close to $\varphi_L$, the energy defines a local shape
	functional
	\[
	T(\varphi):=J(u_{\Omega_\varphi}),
	\]
	where $u_{\Omega_\varphi}$ is the solution of \eqref{problem_Omega} in
	$\Omega_\varphi$ given by the local uniqueness result recalled above, and
	we can consider the pairs $(\Omega_\varphi,u_{\Omega_\varphi})$ which are
	stationary for $T$ under the volume constraint (see
	Section~\ref{section:stability} for the precise definitions). By
	\cite[Proposition 2.6]{AfonsoPacella}, a pair
	$(\Omega_\varphi,u_{\Omega_\varphi})$ is energy-stationary under the
	volume constraint if and only
	if $u_{\Omega_\varphi}$ solves, for some constant $c$, the relative
	overdetermined problem
	\begin{equation}\label{overdetermined}
		\begin{cases}
			-\Delta u=|u|^{p-1}u & \text{in }\Omega_\varphi,\\
			u=0 & \text{on }\Gamma_\varphi,\\
			\dfrac{\partial u}{\partial\nu}=0 & \text{on }\Gamma_{1,\varphi},\\
			\dfrac{\partial u}{\partial\eta}=c & \text{on }\Gamma_\varphi,
		\end{cases}
	\end{equation}
	where $\eta$ is the outer unit normal vector on $\Gamma_{\varphi}$.\\
	Hence, studying energy-stationary pairs amounts to studying the domains
	$\Omega_\varphi$ in which the relative overdetermined problem
	\eqref{overdetermined} admits a solution. In particular, since the
	normal derivative of $u_{p,L}^{(m)}$ is constant on $\Gamma_{\varphi_L}$,
	the pairs $(\Omega_{\varphi_L},\pm u_{p,L}^{(m)})$ are energy-stationary,
	and it is natural to ask whether they locally minimize $T$ among the
	domains $\Omega_\varphi$ with the same volume. In the terminology of
	\cite{AfonsoPacella}, recalled in Section~\ref{section:stability}, this
	amounts to asking whether they are stable or unstable energy-stationary
	pairs.
	
	This question was first addressed in \cite{AfonsoPacella}, where
	stability criteria were obtained for general nonlinearities $f(u)$
	\cite[Theorem 4.11]{AfonsoPacella} and sharp results were proved for the
	torsion problem, $f(u)\equiv1$ \cite[Theorem 1.4]{AfonsoPacella}: in
	this case, the stability of the flat cylinders changes when
	$L^2\lambda_1(\omega)$ crosses the value $\beta\approx1.439$, the unique
	positive solution of $\sqrt\beta\tanh\sqrt\beta=1$, where
	$\lambda_1(\omega)$ is the first nontrivial eigenvalue of
	$-\Delta_{\mathbb R^{N-1}}$ in $\omega$ with homogeneous Neumann boundary
	conditions. In correspondence with this threshold, a bifurcation result has been obtained in \cite{PacellaRuizSicbaldi} (see also \cite{FMW}), while bifurcation for the overdetermined problem was obtained in \cite{LianPacellaSicbaldi}.
	
	For the Lane--Emden equation, the positive one-dimensional solution,
	corresponding to $m=1$, was studied in \cite{DeMPacMaz}, where energy stability for small and large values of $p$ is related to the asymptotic behaviour of the solution and of
	its linearized operator as $p\to1^+$ and $p\to+\infty$. In this paper
	we address the sign-changing case, with an arbitrary number of nodal
	regions; to the best of our knowledge, this is the first stability
	analysis of one-dimensional sign-changing solutions. Our results are
	based on a detailed study of the one-dimensional profiles, mainly
	derived from that of the positive solution, and on the precise
	asymptotic behaviour of the first eigenvalue of the associated
	linearized operator, which instead requires taking into account the
	interaction between the nodal regions. In addition, we obtain
	comparison results valid for the whole range $p>1$ and refine the
	expansions of the amplitude in both regimes, with explicit constants,
	also in the positive case.
	
	Since the even extension of a solution of \eqref{ODE_in_(0,L)} solves
	the same equation in $(-L,L)$ with Dirichlet boundary conditions, and in
	view of the scaling properties of the nonlinearity (see Remark
	\ref{rem:rescaling}), it is enough to study the corresponding Dirichlet
	problem on the fixed interval $I:=(-1,1)$. For $p>1$, let $v_p$ be the
	unique positive solution of
	\begin{equation}\label{positive_sol}
		\begin{cases}
			-v''=v^p & \text{in }I,\\
			v>0 & \text{in }I,\\
			v(\pm1)=0.
		\end{cases}
	\end{equation}
	For $n\geq1$ we set
	\begin{align}
		&r_{n,i}:=-1+\frac{2i}{n},
		\qquad i=0,\ldots,n,\label{ri_si_Ii}\\
		&s_{n,i}:=-1+\frac{2i-1}{n},\qquad
		I_{n,i}:=(r_{n,i-1},r_{n,i})
		=\left(s_{n,i}-\frac1n,s_{n,i}+\frac1n\right),
		\qquad i=1,\ldots,n,\label{si_Ii}
	\end{align}
	and we denote by $v_p^{(n)}$ the solution of
	\begin{equation}\label{problema_in_I}
		\begin{cases}
			-v''=|v|^{p-1}v & \text{in }I,\\
			v(\pm1)=0,
		\end{cases}
	\end{equation}
	having $n$ nodal regions and satisfying $(v_p^{(n)})'(1)<0$.
	In Proposition \ref{prop:existence_uniqueness} we show its existence
	and uniqueness and that it consists of alternating translated and
	rescaled copies of $v_p$; in particular, its zeros are the points
	$r_{n,i}$ and the maximum points of $|v_p^{(n)}|$ are the centers
	$s_{n,i}$ of its nodal intervals $I_{n,i}$. The connection with
	\eqref{ODE_in_(0,L)} is given by even reflection and scaling:
	\[
	u_{p,L}^{(m)}(t)
	=L^{-\frac{2}{p-1}}v_p^{(2m-1)}\!\left(\frac{t}{L}\right),
	\qquad t\in[0,L].
	\]
	Thus a solution with $m$ nodal regions on $(0,L)$ corresponds to a
	Dirichlet solution with $2m-1$ nodal regions on the reflected interval.
	Only odd values $n=2m-1$ enter the stability problem. We nevertheless
	state our results on $v_p^{(n)}$ for every $n\geq1$.
	
	Besides the profiles themselves, a central role is played by the
	first Dirichlet eigenvalue $\widetilde\alpha_1(p,n)$ of the
	linearized operator
	\begin{equation}\label{L_vpn}
		L_{v_p^{(n)}}:=-\frac{d^2}{dt^2}-p|v_p^{(n)}|^{p-1}
		\qquad\text{in }I.
	\end{equation}
	
	We first consider the limit $p\to+\infty$. In addition to extending to
	the nodal setting the asymptotic description obtained in
	\cite{DeMPacMaz}, we determine
	the constant term in the expansion of $\|v_p^{(n)}\|_\infty^{p-1}$,
	which was not available in \cite{DeMPacMaz} even for the positive
	solution and which in turn determines the second-order term in the expansion of $\|\vpn\|_\infty$.
	
	\begin{teorema}[Asymptotic behaviour for $p\to+\infty$]
		\label{Teo_p_to_inf_rescaled}
		Let $n\geq1$ and let $\vpn$ be the solution of
		\eqref{problema_in_I} having $n$ nodal regions. Then, as $p\to+\infty$:
		\begin{enumerate}
			\item\label{item:norm_p-1_inf}
			\[
			\|\vpn\|_\infty^{p-1}=\frac{n^2p}{2}+\frac{n^2(1+4\log2)}{2}
			+O\left(\frac1p\right);
			\]
			\item\label{item:norm_inf}  
			\[
			\|\vpn\|_\infty=1+\frac1p\log\frac{n^2p}{2}+\frac1{2p^2}\log^2\frac{n^2p}{2}
			+\frac1{p^2}\left(\log\frac{n^2p}{2}+1+4\log2\right)
			+O\left(\frac{(\log p)^3}{p^3}\right);
			\]
			in particular,
			$\|\vpn\|_\infty=1+\frac{\log(n^2p/2)}{p}
			+O\left(\frac{(\log p)^2}{p^2}\right)$;
			\item
			\[
			\vpn\longrightarrow v_\infty^{(n)}\quad\text{in }C^0(\bar I)
			\cap C^1_{\rm loc}\bigl(\bar I\setminus\cup_{i=1}^n\{s_{n,i}\}\bigr),
			\]
			where, for any $i\in\{1,\ldots,n\}$,
			\[
			v_\infty^{(n)}(t)=(-1)^{n+i}(1-n|t-s_{n,i}|)\qquad\text{in }\overline{I}_{n,i};
			\]
			\item setting
			\begin{equation}\label{mupn}
				\mupn:=(p\|\vpn\|_\infty^{p-1})^{-\frac12}
			\end{equation}
			and defining, for $i\in\{1,\ldots,n\}$, the rescaling of $\vpn$
			around the maximum point $s_{n,i}$ of $|\vpn|$ in $I_{n,i}$ as
			\begin{equation}\label{utildep_rescaled}
				\vpntildei(s):=p\,\frac{\vpn(\mupn s+s_{n,i})-\vpn(s_{n,i})}{\vpn(s_{n,i})},
				\qquad s\in\left[\tfrac{-1-s_{n,i}}{\mupn},\tfrac{1-s_{n,i}}{\mupn}\right],
			\end{equation}
			we have, for every $i\in\{1,\ldots,n\}$,
			\[
			\vpntildei\longrightarrow W\quad\text{in }C^1_{\rm loc}(\mathbb R),
			\]
			where
			\begin{equation}\label{W}
				W(s)=\log\frac{4e^{\sqrt{2}s}}{(1+e^{\sqrt{2}s})^2},\quad s\in\mathbb{R},
			\end{equation}
			is the solution to the limit problem
			\begin{equation}\notag
				\begin{cases}
					-W''=e^W & \text{in }\mathbb{R},\\
					W'(0)=0,\\
					W(0)=0;
				\end{cases}
			\end{equation}
			\item
			\[
			\tilde\alpha_1(p,n)=-\frac{n^2p^2}{4}(1+o_p(1)),
			\]
			where $\widetilde\alpha_1(p,n)$ is the first eigenvalue of $\eqref{L_vpn}$.
		\end{enumerate}
	\end{teorema}
	
	We next consider the opposite asymptotic regime $p\to1^+$. 
	
	\begin{teorema}[Asymptotic behaviour for $p\to1^+$]
		\label{Teo_p_to_1_rescaled}
		Let $n\geq1$ and let $\vpn$ be the solution of
		\eqref{problema_in_I} having $n$ nodal regions. Then, as $p\to1^+$:
		\begin{enumerate}
			\item
			\[
			\|\vpn\|_\infty^{p-1}
			=n^2\frac{\pi^2}{4}
			+n^2\frac{\pi^2}{4}\left(\log2-\frac12\right)(p-1)+O\bigl((p-1)^2\bigr);
			\]
			\item
			\[
			\frac{\vpn}{\|\vpn\|_\infty}\longrightarrow\varphi_n
			\qquad\text{in }C^1(\bar I),
			\]
			where $\varphi_n$ is the eigenfunction corresponding to
			$\nu_n=n^2\frac{\pi^2}{4}$, the $n$-th Dirichlet eigenvalue of
			$-\frac{d^2}{dt^2}$ in $I$, normalized by $\|\varphi_n\|_\infty=1$ and
			$\varphi_n'(1)<0$, namely
			\[
			\varphi_n(t)=
			\begin{cases}
				(-1)^{\frac{n-1}{2}}\cos\left(\dfrac{n\pi}{2}t\right), & \text{if $n$ is odd},\\[2mm]
				(-1)^{\frac n2+1}\sin\left(\dfrac{n\pi}{2}t\right), & \text{if $n$ is even};
			\end{cases}
			\]
			\item
			\[
			p|\vpn|^{p-1}\longrightarrow n^2\frac{\pi^2}{4}
			\qquad\text{in }C^1_{\rm loc}\bigl(I\setminus\{r_{n,1},\ldots,r_{n,n-1}\}\bigr),
			\]
			where the set of interior zeros is empty if $n=1$;
			\item
			\begin{equation}\label{converg_alpha1p}
				\tilde\alpha_1(p,n)\longrightarrow(1-n^2)\frac{\pi^2}{4},
			\end{equation}
			where $\widetilde\alpha_1(p,n)$ is the first eigenvalue of $\eqref{L_vpn}$.
		\end{enumerate}
	\end{teorema}
	
	Besides providing the asymptotic information needed for the stability
	problem, we also obtain a global description of the dependence of the
	one-dimensional solutions and of the associated linearized potentials
	on the exponent $p$, which is of independent interest. In addition to
	an exact formula for the $\|\vpn\|_\infty$, we prove a strict pointwise
	comparison between the absolute values of the solutions for different
	exponents, and we describe precisely the intersections of the
	associated linearized potentials. Set
	\begin{equation}\label{J_p}
		J_p:=\int_0^1\frac{ds}{\sqrt{1-s^{p+1}}}.
	\end{equation}
	
	\begin{teorema}\label{teo:globale_any_n}
		Let $n\geq1$ and, for $p>1$, let $\vpn$ be the solution of
		\eqref{problema_in_I} having $n$ nodal regions. Then:
		\begin{enumerate}
			\item for every $p>1$,
			\[
			\|\vpn\|_\infty^{p-1}
			=
			\frac{n^2(p+1)}{2}J_p^2;
			\]
			
			\item if $p_2>p_1>1$, then, for every
			$i\in\{1,\ldots,n\}$,
			\[
			|v_{p_2}^{(n)}(t)|
			<
			|v_{p_1}^{(n)}(t)|
			\qquad
			\text{for every }t\in I_{n,i}.
			\]
			In particular, the map
			$p\mapsto\|\vpn\|_\infty$ is strictly decreasing in
			$(1,+\infty)$;
			
			\item if $p_2>p_1>1$, then there exists a unique
			$\tau=\tau(p_1,p_2)\in(0,1)$ such that, for every
			$i\in\{1,\ldots,n\}$,
			\begin{align*}
				p_2|v_{p_2}^{(n)}(t)|^{p_2-1}
				&>
				p_1|v_{p_1}^{(n)}(t)|^{p_1-1}
				&&\text{if }
				|t-s_{n,i}|<\frac{\tau}{n},\\
				p_2|v_{p_2}^{(n)}(t)|^{p_2-1}
				&<
				p_1|v_{p_1}^{(n)}(t)|^{p_1-1}
				&&\text{if }
				\frac{\tau}{n}<|t-s_{n,i}|<\frac1n.
			\end{align*}
			Consequently, the two potentials have exactly two intersection
			points in each $I_{n,i}$, namely
			$s_{n,i}-\frac{\tau}{n}$ and $s_{n,i}+\frac{\tau}{n}$.
		\end{enumerate}
	\end{teorema}
	
	The exact formula in (1), combined with the representation of $J_p$ in
	terms of the Beta function, is the common source of the amplitude
	expansions in Theorem \ref{Teo_p_to_inf_rescaled} (1)--(2) and Theorem
	\ref{Teo_p_to_1_rescaled} (1). 
	To the best of our knowledge, items (2) and (3) are new even for the
	positive solution $v_p$. Moreover, the monotonicity in (2) is specific
	to the one-dimensional setting. Indeed, let $w_{p,N}$ denote the positive (hence radial) solution of the Lane--Emden
	problem in the unit ball $B_1\subset\mathbb R^N$ with homogeneous
	Dirichlet boundary conditions. Then $\|w_{p,N}\|_\infty\to+\infty$ as
	$p\to1^+$, as can be easily seen by testing the equation against a
	positive first Dirichlet eigenfunction of $-\Delta$ in $B_1$. On the other hand, for $N=2$ the $L^\infty$-norm converges to
	$\sqrt e$ from below as $p\to+\infty$ (see \cite{GrossiIanniLuoYan}), while
	for $N\geq3$ it diverges as $p\to\bigl(\frac{N+2}{N-2}\bigr)^-$
	(see \cite{AtkinsonPeletier}). Hence, for every $N\geq2$, the map
	$p\mapsto\|w_{p,N}\|_\infty$ is not monotone.
	
	We now turn to energy stability in the flat cylinder. As it will be
	clear from the statements below, our results do not depend on the shape
	of the base $\omega$, which enters only through its first nontrivial
	Neumann eigenvalue $\lambda_1(\omega)$, namely the first positive
	eigenvalue of
	\[
	\begin{cases}
		-\Delta_{\mathbb R^{N-1}}\psi=\lambda\psi & \text{in }\omega,\\
		\dfrac{\partial\psi}{\partial\nu}=0 & \text{on }\partial\omega.
	\end{cases}
	\]
	In particular, no symmetry assumption on $\omega$ is required. After
	rescaling the height to one, the relevant geometric parameter is
	$L^2\lambda_1(\omega)$. As explained in Section \ref{section:stability},
	a sufficient condition for the nondegeneracy of $u_{p,L}^{(m)}$ is
	\begin{equation}\label{nondeg_cond}
		L^2\lambda_1(\omega)>-\widetilde\alpha_1(p,2m-1).
	\end{equation}
	Under this condition, the constrained second variation can be analysed
	through an auxiliary one-dimensional linear problem, using the
	characterization in \cite[Theorem 4.11]{AfonsoPacella}, which we recall
	in Theorem \ref{teorema}. The asymptotic information
	above then yields explicit stability and instability regimes. For large
	exponents, we obtain the following sufficient condition for stability.
	
	\begin{teorema}[Stability for $p\to+\infty$]\label{Teo_stability_p_infty}
		Let $m\geq1$ and $\gamma>m-\frac12$. There exists $p_0=p_0(m,\gamma)>1$
		such that the following holds. Let $\omega\subset\mathbb R^{N-1}$ be a
		smooth bounded domain, let $p\geq p_0$ and let
		\[
		L>\frac{\gamma p}{\sqrt{\lambda_1(\omega)}}.
		\]
		Then the solution $u_{p,L}^{(m)}$ defined in \eqref{write_sol} is
		nondegenerate and the pairs $(\Omega_{\varphi_L},\pm u_{p,L}^{(m)})$ are
		stable energy-stationary pairs in the sense of
		Definition~\ref{def:critical_energy_pair_stable}.
	\end{teorema}
	
	For exponents close to one, both stability and instability occur,
	depending on the height of the cylinder.
	
	\begin{teorema}[Stability and instability for $p\to1^+$]\label{Teo_stability_p_1}
		Let $m\geq1$, let $\omega\subset\mathbb R^{N-1}$ be a smooth bounded
		domain, let $L>0$ and let $u_{p,L}^{(m)}$ be defined in \eqref{write_sol}.
		\begin{enumerate}
			\item If
			\[
			L>\frac{\pi(2m-1)}{2\sqrt{\lambda_1(\omega)}},
			\]
			there exists $p_0>1$, depending only on $m$ and $L^2\lambda_1(\omega)$,
			such that, for every $p\in(1,p_0)$, $u_{p,L}^{(m)}$ is nondegenerate and
			the pairs $(\Omega_{\varphi_L},\pm u_{p,L}^{(m)})$ are stable
			energy-stationary pairs.
			\item If
			\[
			\frac{\pi\sqrt{m(m-1)}}{\sqrt{\lambda_1(\omega)}}
			<L<\frac{\pi(2m-1)}{2\sqrt{\lambda_1(\omega)}},
			\]
			there exists $p_0>1$, depending only on $m$ and $L^2\lambda_1(\omega)$,
			such that, for every $p\in(1,p_0)$, $u_{p,L}^{(m)}$ is nondegenerate and
			the pairs $(\Omega_{\varphi_L},\pm u_{p,L}^{(m)})$ are unstable
			energy-stationary pairs.
		\end{enumerate}
		Stability and instability are understood in the sense of Definitions
		\ref{def:critical_energy_pair_stable} and
		\ref{def:critical_energy_pair_unstable}, respectively.
	\end{teorema}
	
	The two height scales in the latter theorem have distinct roles.
	The lower endpoint of the instability interval comes from the
	limiting sufficient condition for nondegeneracy, since
	$\widetilde\alpha_1(p,2m-1)\to -m(m-1)\pi^2$ as $p\to1^+$.
	Within this regime, the stability and instability conditions are
	separated by the limiting height
	$\pi(m-\tfrac12)/\sqrt{\lambda_1(\omega)}$.
	For large $p$, the sufficient stability condition instead involves
	heights growing linearly with the exponent. In both regimes, the
	number of nodal regions appears explicitly in the geometric
	conditions for stability. Notice that in the case $m=1$ these theorems reduce to \cite[Theorems 1.3 and 1.4]{DeMPacMaz}.
	
	Theorems \ref{Teo_stability_p_infty} and \ref{Teo_stability_p_1}
	describe the stability of the flat cylinders only in the asymptotic
	regimes $p\to+\infty$ and $p\to1^+$. It would be interesting to
	determine whether, for every $p>1$ and every $m\geq1$, there exists a
	threshold for the geometric parameter $L^2\lambda_1(\omega)$ separating
	stable and unstable energy-stationary pairs
	$(\Omega_{\varphi_L},\pm u_{p,L}^{(m)})$, and to describe its dependence
	on $p$; by Theorems \ref{Teo_stability_p_1} and
	\ref{Teo_stability_p_infty}, such a threshold should converge to
	$(m-\frac12)^2\pi^2$ as $p\to1^+$ and grow at most quadratically in $p$
	as $p\to+\infty$. Moreover, in analogy with
	\cite{PacellaRuizSicbaldi, LianPacellaSicbaldi}, one may expect that,
	at this threshold, a branch of nontrivial domains $\Omega_\varphi$
	admitting solutions of \eqref{overdetermined} bifurcates from the flat
	cylinders. We plan to address these questions in a forthcoming work,
	possibly also for more general nonlinearities.
	
	A bifurcation result of a different nature was obtained in
	\cite{Afonso1} for a class of semilinear problems in flat cylinders,
	which includes the Lane--Emden one. There the domain is not deformed:
	the flat cylinder is only rescaled, and bifurcation occurs at the values
	of the parameter for which the one-dimensional solution becomes
	degenerate, under a simplicity assumption which holds for generic bases
	$\omega$. The bifurcating solutions are solutions of
	\eqref{problem_Omega} in flat cylinders which are not one-dimensional
	and have the same number of nodal regions as the one-dimensional
	solution; in particular, they are not solutions of
	\eqref{overdetermined}.
	
	The paper is organized as follows. In Section~\ref{section:positive} we
	study the positive solution $v_p$, establishing the global comparison
	results and refining its asymptotic expansions as $p\to+\infty$ and as
	$p\to1^+$. In Section~\ref{section:nodal} we construct the nodal
	solutions of \eqref{problema_in_I}, prove Theorems
	\ref{Teo_p_to_inf_rescaled}, \ref{Teo_p_to_1_rescaled} and
	\ref{teo:globale_any_n}, and analyse the first eigenvalue of the
	corresponding linearized operator. Section~\ref{section:stability}
	contains the variational framework and the proofs of Theorems
	\ref{Teo_stability_p_infty} and \ref{Teo_stability_p_1}.

	\section{Preliminaries and results on the positive solution}\label{section:positive}
	
	In this section we collect some properties of the positive solution $\vp$ introduced in \eqref{positive_sol}. It is well known that $v_p$ is even, decreasing in $(0,1)$ and thus $v_p'(0)=0$ and $\|v_p\|_\infty=v_p(0)$.
	\subsection{Dependence on the exponent}
	\begin{proposizione}\label{prop:monotonia_L_infty}
		Let $\vp$ be the positive solution of \eqref{positive_sol}, then we have
		\begin{enumerate}
			\item[(i)]
			\[
			\|v_p\|_\infty^{p-1}
			=
			\frac{p+1}{2}J_p^2\qquad \text{for any $p>1$},
			\]
			where $J_p$ is defined in \eqref{J_p};
			\item[(ii)] the map
			\[
			p\longmapsto \|v_p\|_\infty
			\]
			is strictly decreasing in $(1,+\infty)$.
		\end{enumerate}
	\end{proposizione}
	
	\begin{proof}
		Set
		\[
		M_p:=\|v_p\|_\infty=v_p(0).
		\]
		Multiplying the equation in \eqref{positive_sol} by $v_p'$ and integrating from $0$ to $x$, we obtain
		\[
		\frac12 (v_p'(x))^2
		=
		\frac{1}{p+1}
		\left(M_p^{p+1}-v_p(x)^{p+1}\right).
		\]
		Since $v_p'<0$ in $(0,1)$, it follows that for $x\in(0,1)$
		\[
		-v_p'(x)
		=
		\sqrt{\frac{2}{p+1}}\,
		\sqrt{M_p^{p+1}-v_p^{p+1}(x)}.
		\]
		Separating the variables and integrating,
		we get
		\[
		1
		=
		\sqrt{\frac{p+1}{2}}
		\int_0^{M_p}
		\frac{dv}{\sqrt{M_p^{p+1}-v^{p+1}}}.
		\]
		With the change of variables $v=M_p s$, this becomes
		\[
		1
		=
		\sqrt{\frac{p+1}{2}}\,
		M_p^{-\frac{p-1}{2}}
		\int_0^1\frac{ds}{\sqrt{1-s^{p+1}}}
		=
		\sqrt{\frac{p+1}{2}}\,
		M_p^{-\frac{p-1}{2}}J_p.
		\]
		Therefore
		\[
		M_p^{p-1}
		=
		\frac{p+1}{2}J_p^2.
		\]
		This concludes the proof of (i).
		
		We now prove (ii). By (i), we have
		\begin{equation}\label{da stimare}
			\log M_p
			=
			\frac{1}{p-1}\log\frac{p+1}{2}
			+
			\frac{2}{p-1}\log J_p.
		\end{equation}
		The first term in the right-hand side is strictly decreasing in $p$.
		Indeed,
		\begin{equation}\label{primo_termine}
			\frac{d}{dp}
			\left(
			\frac{1}{p-1}\log\frac{p+1}{2}
			\right)
			=
			\frac{1}{(p-1)^2}
			\left(
			\frac{p-1}{p+1}
			-
			\log\frac{p+1}{2}
			\right)<0,
		\end{equation}
		where the last inequality follows from $x\log x>x-1$, for $x>1$, applied with $x=\frac{p+1}{2}$. 
		
		Moreover $J_p>1$ and $J_p$ is strictly decreasing in $p$, as follows
		immediately from its definition. Hence
		\begin{equation}\label{secondo_termine}
			p\longmapsto \frac{2}{p-1}\log J_p
		\end{equation}
		is strictly decreasing as well. Combining \eqref{da stimare}, \eqref{primo_termine} and \eqref{secondo_termine} we deduce that $p\mapsto\log M_p$ is strictly decreasing in $p$, and thus the thesis follows.
	\end{proof}

	\begin{proposizione}\label{prop:ordine_potenziale}
		For $p>1$, let $\vp$ be the positive solution of \eqref{positive_sol}.\\ Given $p_2>p_1>1$, there exists a unique
		$\tau=\tau(p_1,p_2)\in(0,1)$ such that
		\[
		p_2 v_{p_2}^{p_2-1}(t)>p_1 v_{p_1}^{p_1-1}(t)
		\quad\text{if }|t|\in[0,\tau),
		\qquad
		p_2 v_{p_2}^{p_2-1}(t)<p_1 v_{p_1}^{p_1-1}(t)
		\quad\text{if }|t|\in(\tau,1).
		\]
		Furthermore,
		\[
		v_{p_2}<v_{p_1}
		\qquad\text{in }I.
		\]
	\end{proposizione}
	
	\begin{proof} Clearly by the evenness of the solutions $v_p$, it suffices to prove both assertions in $[0,1)$.\\
		Set
		\[
		M_p:=\|v_p\|_\infty=v_p(0)\qquad\text{and}\qquad
		y_p:=\frac{v_p}{M_p}.
		\]
		Reasoning as in the proof of Proposition \ref{prop:monotonia_L_infty}, we have, for $x\in[0,1)$
		\[
		-v_p'(x)
		=
		\sqrt{\frac{2}{p+1}}\,
		\sqrt{M_p^{p+1}-v_p^{p+1}(x)}.
		\]
		Therefore
		\[
		-y_p'(x)
		=
		\sqrt{\frac{2}{p+1}}\,
		M_p^{\frac{p-1}{2}}
		\sqrt{1-y^{p+1}_p(x)},
		\]
		moreover by Proposition \ref{prop:monotonia_L_infty} (i),
		\[
		\frac1{J_p}\frac{-y'_p(x)}{\sqrt{1-y_p^{p+1}(x)}}
		=1,
		\]
		where $J_p$ is defined in \eqref{J_p}.
		Consequently, for $t\in(0,1)$, integrating in $[0,t]$ we get
		\begin{equation}\label{separazione_variabili}
			t
			=
			\frac1{J_p}
			\int_{y_p(t)}^1
			\frac{ds}{\sqrt{1-s^{p+1}}},
		\end{equation}
		where $J_p$ is defined in \eqref{J_p}.\\
		Let
		\[
		\rho_p(s):=
		\frac1{J_p\sqrt{1-s^{p+1}}},
		\qquad s\in[0,1).
		\]
		By definition $\rho_p(0)=\frac{1}{J_p}$, then, being $J_p$ strictly monotone decreasing in $p$, given $p_2>p_1>1$
		\begin{equation}\label{in_0}
			\rho_{p_2}(0)>\rho_{p_1}(0). 
		\end{equation}
		Moreover
		\[
		\frac{\rho_{p_2}(s)}{\rho_{p_1}(s)}
		=
		\frac{J_{p_1}}{J_{p_2}}
		\sqrt{\frac{1-s^{p_1+1}}{1-s^{p_2+1}}}
		\]
		is strictly decreasing in $s$, for $s\in(0,1)$, being
		\begin{equation}\label{derivative}
			\frac{d}{ds}\left(\frac{1-s^{p_1+1}}{1-s^{p_2+1}}\right)
			=
			\frac{s^{p_1}}
			{(1-s^{p_2+1})^2}
			\left[
			-(p_1+1)
			+(p_2+1)s^{p_2-p_1}
			-(p_2-p_1)s^{p_2+1}
			\right]<0.
		\end{equation}
		Indeed, the expression in square brackets vanishes at $s=1$ and its
		derivative verifies
		\[
		(p_2+1)(p_2-p_1)s^{p_2-p_1-1}
		\left(1-s^{p_1+1}\right)>0
		\qquad\text{in }(0,1),
		\]
		hence it is strictly negative in $(0,1)$.
		
		By definition of $J_p$, see $\eqref{J_p}$,
		\begin{equation}\label{integral}
			1=\int_0^1 \rho_p(s)\,ds\qquad\text{for any $p>1$},
		\end{equation}
		so by \eqref{in_0} and \eqref{derivative} we derive that $\frac{\rho_{p_2}}{\rho_{p_1}}$ crosses the value $1$ exactly once and it is smaller than $1$ at $s=1$.\\
		These latter facts, combined with \eqref{integral}, imply that
		\[
		\int_y^1 \rho_{p_2}(s)\,ds<\int_y^1  \rho_{p_1}(s)\,ds\qquad\text{for any $y\in(0,1)$.}
		\]
		Together with \eqref{separazione_variabili}, this gives 
		\begin{equation}\label{yp_ordine}
			y_{p_2}(t)<y_{p_1}(t)
			\qquad\text{for any }t\in(0,1),
		\end{equation}
		indeed 
		\[
		\int_{y_{p_1}(t)}^1 \rho_{p_1}(s)ds
		=\int_{y_{p_2}(t)}^1 \rho_{p_2}(s)ds
		<\int_{y_{p_2}(t)}^1 \rho_{p_1}(s)ds
		\]
		and then \eqref{yp_ordine} follows from the positivity of $\rho_{p_1}$.
		
		Next we observe that the quantity $(p-1)J_p$ is strictly increasing in $p$, as it can be easily seen via the change of variables $t=s^{p+1}$, indeed
		\[
		(p-1)J_p=\int_0^1\left(1-\frac2{p+1}\right)\frac{t^{\frac1{p+1}}}{t\sqrt{1-t}}dt.
		\]
		
		Moreover, for $0<y<1$, the function
		\[
		F_p(y):=\frac{\sqrt{1-y^{p+1}}}{y}
		\]
		is strictly decreasing in $y$ and strictly increasing in $p$.
		Hence, by \eqref{yp_ordine},
		$F_{p_2}(y_{p_2}(x))
		>F_{p_1}(y_{p_1}(x))$.\\
		Combining the monotonicity properties established above, we obtain 
		\[
		(p_2-1)J_{p_2}F_{p_2}(y_{p_2}(x))
		>
		(p_1-1)J_{p_1}F_{p_1}(y_{p_1}(x)),
		\]
		and therefore
		\[
		\frac{d}{dx}\left(\log \left(\frac{p_2v_{p_2}^{p_2-1}(x)}{p_1v_{p_1}^{p_1-1}(x)}\right)\right)
		=(p_1-1)J_{p_1}
		\frac{\sqrt{1-y_{p_1}(x)^{p_1+1}}}{y_{p_1}(x)}-(p_2-1)J_{p_2}
		\frac{\sqrt{1-y_{p_2}(x)^{p_2+1}}}{y_{p_2}(x)}<0.
		\]
		In conclusion we have that
		\begin{equation}\label{decreasing}
			\frac{p_2v_{p_2}^{p_2-1}}{p_1v_{p_1}^{p_1-1}}
			\quad\text{is strictly decreasing in }(0,1).
		\end{equation}
		
		By Proposition \ref{prop:monotonia_L_infty} (i),
		\[
		pM_p^{p-1}
		=
		\frac{p(p+1)}2J_p^2,
		\]
		then, using again the above change of variables, it can be easily seen that $pM_p^{p-1}$ is strictly increasing in $p$, being
		\[
		pM_p^{p-1}
		=\frac12\left(1-\frac1{p+1}\right)
		\left(\int_0^1\frac{t^{\frac1{p+1}}}{t\sqrt{1-t}}dt\right)^2.
		\]
		In particular
		\begin{equation}\label{at_0}
			p_2v_{p_2}^{p_2-1}(0)>p_1v_{p_1}^{p_1-1}(0).
		\end{equation}
		On the other hand, since $v_{p_i}'(1)<0$, for $i=1,2$,
		$v_{p_i}(x)=-v_{p_i}'(1)(1-x)+o(1-x)$ as $x\to1^-$, and hence as $x\to1^-$
		\begin{equation}\label{at_1}
			\frac{p_2v_{p_2}^{p_2-1}(x)}{p_1v_{p_1}^{p_1-1}(x)}=C(1-x)^{p_2-p_1}(1+o(1))
			\longrightarrow0.
		\end{equation}
		By \eqref{decreasing}, \eqref{at_1} and \eqref{at_0} there exists
		a unique $\tau\in(0,1)$ such that $
		p_2v_{p_2}^{p_2-1}>p_1v_{p_1}^{p_1-1}$ in $[0,\tau)$, and $p_2v_{p_2}^{p_2-1}<p_1v_{p_1}^{p_1-1}$
		in $(\tau,1)$.
		
		Finally, from Proposition \ref{prop:monotonia_L_infty} \emph{(ii)} and \eqref{yp_ordine} we get
		\[
		v_{p_2}=M_{p_2}y_{p_2}<M_{p_1}y_{p_1}=v_{p_1}\qquad\text{in }[0,1)
		\]
		which concludes the proof.
	\end{proof}

	\subsection{Asymptotic behaviour}
	Let $\tilde\alpha_1(p)$ be the first eigenvalue of the linearized operator at $\vp$ with Dirichlet boundary conditions in $I$, defined as
	\begin{equation}\label{L_vp}
		L_{\vp}=-\frac{d^2}{dt^2}-p\vp^{p-1}.
	\end{equation}
	We recall some results obtained in \cite{DeMPacMaz}, concerning the asymptotic behaviour of $v_p$ and $\tilde\alpha_1(p)$ as $p\to+\infty$ and as $p\to1$ and we refine the asymptotic expansions of $\|v_p\|_\infty$ as $p\to+\infty$ and of $\|v_p\|_\infty^{p-1}$ as $p\to1^+$.
	
	\begin{proposizione}[{\cite[Theorem 1.1 and Lemma 2.9]{DeMPacMaz}}]\label{Teo_p_to_inf}
		Let $\vp$ be the positive solution of \eqref{positive_sol}, then we have
		\begin{enumerate}
			\item \begin{equation*}
				\lim_{p\to +\infty}\Vert \vp \Vert_{\infty}=1;
			\end{equation*}
			\item \begin{equation*}
				\vp\rightarrow 2G(\cdot,0)=1-|\cdot| \qquad \mbox{ in } C^0([-1,1])\cap C^1_{loc}([-1,1]\setminus\lbrace 0 \rbrace),
			\end{equation*}
			where $G$ is the Green's function of $-\frac{d^2}{d t^2}$ with Dirichlet boundary conditions defined as
			\begin{equation}\label{Green}
				G(t,\tau)=\begin{cases}
					-\frac{(\tau+1)(t-1)}{2} \quad &\mbox{ if } t>\tau\\
					\frac{(t+1)(1-\tau)}{2}  \quad &\mbox{ if } t\leq \tau;
				\end{cases}
			\end{equation}
			\item \label{as_behav_norm}
			\begin{equation*}
				\Vert \vp \Vert_{\infty}^{p+1}=\frac{p}{2}(1+o_p(1)) \qquad \mbox{as $p\to+\infty$};
			\end{equation*}
			\item\label{mu_p_piu} setting
			\begin{equation}\label{mu_p}
				\mu_p:=(p\|\vp\|_\infty^{p-1})^{-\frac12}
			\end{equation}
			and
			\begin{equation}\label{utildep}
				\tilde{v}_p(s)=p\frac{\vp(\mu_ps)-\|\vp\|_{\infty}}{\|\vp\|_{\infty}}, \qquad \text{for $s \in [-\tfrac{1}{\mu_p}, \tfrac{1}{\mu_p}]$,}
			\end{equation}
			then
			\begin{equation*}
				\tilde{v}_p \rightarrow W \quad \mbox{ in } C^1_{loc}(\mathbb{R})\qquad \mbox{as $p\to+\infty$}
			\end{equation*}
			where 
			\begin{equation}\label{W_rescaled}
				W(s)=\log\frac{4e^{\sqrt{2}s}}{(1+e^{\sqrt{2}s})^2},
				\quad s\in \mathbb{R},
			\end{equation}
			is the solution of the limit problem 
			\begin{equation}\label{eq_W}
				\begin{cases}
					-W''=e^W \quad \mbox{ in } \mathbb{R}\\
					W'(0)=0\\
					W(0)=0;
				\end{cases}
			\end{equation}
			\item\label{item:global_estimate_vtildep} there exist $C_1,\,C_2>0$ such that, for $p$ large enough,
			\begin{equation*}
				\tilde{v}_p \leq C_1 W+C_2 \quad \mbox{ in } \left [-\frac{1}{\mu _p}, \frac{1}{\mu _p} \right];
			\end{equation*}
			\item\label{item:asbehavalpha1p}
			\begin{equation*}
				\tilde\alpha_1(p)=-\frac{1}{4}p^2(1+o_p(1))\qquad \mbox{as $p\to+\infty$},
			\end{equation*}
			where $\tilde\alpha_1(p)$ is the first eigenvalue of \eqref{L_vp} with Dirichlet boundary conditions in $I$.
		\end{enumerate}
	\end{proposizione}
	
	The next result identifies the constant term in the expansion of $\|v_p\|_\infty^{p-1}$, improving \eqref{as_behav_norm} of Proposition \ref{Teo_p_to_inf}, and derives from it the expansion of $\|v_p\|_\infty$ up to second order.
	
	\begin{proposizione}\label{prop:new_as_behav_norm}
		Let $\vp$ be the positive solution of \eqref{positive_sol}. Then, as $p\to+\infty$,
		\begin{equation}\label{migliore_esp_p_to_infty}
			\|\vp\|_\infty^{p-1}=\frac{p}{2}+\frac{1+4\log2}{2}+O\left(\frac1p\right)
		\end{equation}
		and
		\begin{equation}\label{vp_norm_expansion}
			\|\vp\|_\infty=1+\frac1p\log\frac p2+\frac1{2p^2}\log^2\frac p2+\frac1{p^2}\left(\log\frac p2+1+4\log2\right)
			+O\left(\frac{(\log p)^3}{p^3}\right).
		\end{equation}
		In particular,
		$\|\vp\|_\infty=1+\frac{\log(p/2)}{p}+O\left(\frac{(\log p)^2}{p^2}\right)$.
	\end{proposizione}
	
	\begin{proof}
		Set $\varepsilon_p:=\frac1{p+1}$. With the change of variables
		$s=t^{\varepsilon_p}$ in \eqref{J_p} we obtain
		\[
		J_p=\varepsilon_p\int_0^1\frac{t^{\varepsilon_p}}{t\sqrt{1-t}}\,dt
		=\varepsilon_p B\left(\varepsilon_p,\frac12\right)
		=\frac{\Gamma(1+\varepsilon_p)\Gamma(\frac12)}{\Gamma(\frac12+\varepsilon_p)},
		\]
		where $B$ and $\Gamma$ denote the Beta and Gamma functions, and we used
		$B(a,b)=\frac{\Gamma(a)\Gamma(b)}{\Gamma(a+b)}$ and
		$\varepsilon\Gamma(\varepsilon)=\Gamma(1+\varepsilon)$. The function
		\[
		g(\varepsilon):=\log\Gamma(1+\varepsilon)-\log\Gamma\left(\tfrac12+\varepsilon\right)
		+\log\Gamma\left(\tfrac12\right)
		\]
		is smooth in a neighbourhood of $\varepsilon=0$, with $g(0)=0$ and
		$g'(0)=(\log\Gamma)'(1)-(\log\Gamma)'(\frac12)=2\log2$ (see
		\cite[Chapter 6]{AbramowitzStegun}). Since
		$\log J_p=g(\varepsilon_p)$, a Taylor expansion gives
		\begin{equation}\label{logJp_expansion}
			2\log J_p=\frac{4\log2}{p+1}+O\left(\frac1{p^2}\right)
			=\frac{4\log2}{p}+O\left(\frac1{p^2}\right),
		\end{equation}
		and therefore
		\begin{equation}\label{Jp_first_expansion}
			J_p^2=1+\frac{4\log2}{p+1}+O\left(\frac1{p^2}\right).
		\end{equation}
		By Proposition \ref{prop:monotonia_L_infty} (i) and \eqref{Jp_first_expansion},
		\[
		\|\vp\|_\infty^{p-1}=\frac{p+1}{2}J_p^2=\frac{p+1}{2}+2\log2+O\left(\frac1p\right),
		\]
		which is \eqref{migliore_esp_p_to_infty}.
		
		Next, taking logarithms in Proposition \ref{prop:monotonia_L_infty} (i),
		\begin{equation}\label{log_norm_vp}
			\log\|\vp\|_\infty=\frac1{p-1}\left(\log\frac{p+1}{2}+2\log J_p\right).
		\end{equation}
		Since $\log\frac{p+1}{2}=\log\frac p2+\log\left(1+\frac1p\right)
		=\log\frac p2+\frac1p+O\left(\frac1{p^2}\right)$, by \eqref{logJp_expansion}
		the bracket in \eqref{log_norm_vp} equals
		\begin{equation}\label{bracket_expansion}
			\log\frac p2+\frac{1+4\log2}{p}+O\left(\frac1{p^2}\right).
		\end{equation}
		Using $\frac1{p-1}=\frac1p+\frac1{p^2}+O\left(\frac1{p^3}\right)$ and
		$\log\frac p2=O(\log p)$, we obtain
		\[
		x_p:=\log\|\vp\|_\infty=\frac1p\log\frac p2+\frac1{p^2}\left(\log\frac p2+1+4\log2\right)
		+O\left(\frac{\log p}{p^3}\right).
		\]
		In particular $x_p=O\left(\frac{\log p}{p}\right)$, so that
		$x_p^2=\frac1{p^2}\log^2\frac p2+O\left(\frac{(\log p)^2}{p^3}\right)$ and
		$x_p^3=O\left(\frac{(\log p)^3}{p^3}\right)$. Hence
		$\|\vp\|_\infty=e^{x_p}=1+x_p+\frac{x_p^2}{2}+O(x_p^3)$ gives
		\eqref{vp_norm_expansion}. The last assertion follows since
		$\frac1{2p^2}\log^2\frac p2+\frac1{p^2}\left(\log\frac p2+1+4\log2\right)
		=O\left(\frac{(\log p)^2}{p^2}\right)$.
	\end{proof}
	
	\begin{proposizione}[{\cite[Theorem 1.2]{DeMPacMaz}}]\label{Teo_p_to_1}
		Let $\vp$ be the positive solution of \eqref{positive_sol} and let $\varphi_1$ be the first eigenfunction of $-\frac{d^2}{d t^2}$ in $I$ with Dirichlet boundary condition, namely $\varphi_1(t)=\cos\left(\frac{\pi}{2}t\right)$ for $t\in\bar I$, then we have\vspace{0.2cm}
		
		\begin{enumerate}
			\item \begin{equation*}
				\Vert \vp \Vert_{\infty}^{p-1}=\frac{\pi^2}{4}+\frac{\pi^2}{4}\tilde c (p-1)+o(p-1)\qquad\text{as $p \to 1$},
			\end{equation*} 
			where
			\begin{equation*}
				\tilde{c}=\frac{\int_{-1}^{1}\varphi_1^2(t)|\log\left(\varphi_1(t)\right)|\,dt}{\int_{-1}^{1}\varphi_1^2(t)\,dt}>0;
			\end{equation*}
			\item\begin{equation*}
				\lim_{p \to 1}\frac{\vp}{\|\vp\|}_{\infty}=\varphi_1\qquad \text{in $C^1(\bar I)$};
			\end{equation*}
			\item \begin{equation*}
				p\vp^{p-1}\to\frac{\pi^2}{4}\qquad\text{in $C^1_{loc}(I)$\qquad\text{as $p \to 1$};}
			\end{equation*}
			\item \begin{equation*}
				\lim_{p \to 1}\tilde\alpha_1(p)=0,
			\end{equation*}
			where $\tilde\alpha_1(p)$ is the first eigenvalue of \eqref{L_vp}.
		\end{enumerate}
	\end{proposizione}
	
	The exact formula in Proposition \ref{prop:monotonia_L_infty} (i) also gives the expansion of $\|\vp\|_\infty^{p-1}$ as $p\to1^+$, with an explicit constant.
	
	\begin{proposizione}\label{prop:explicit_c}
		Let $\vp$ be the positive solution of \eqref{positive_sol}. Then, as $p\to1^+$,
		\begin{equation}\label{esp_p_to_1}
			\|\vp\|_\infty^{p-1}=\frac{\pi^2}{4}+\frac{\pi^2}{4}\left(\log2-\frac12\right)(p-1)
			+O\bigl((p-1)^2\bigr).
		\end{equation}
		In particular, the constant $\tilde c$ in Proposition \ref{Teo_p_to_1} (1)
		is $\tilde c=\log2-\frac12$.
	\end{proposizione}
	
	\begin{proof}
		As in the proof of Proposition \ref{prop:new_as_behav_norm}, setting
		$\varepsilon_p:=\frac1{p+1}$ we have $J_p=\varepsilon_pB(\varepsilon_p,\frac12)$.
		Hence, by Proposition \ref{prop:monotonia_L_infty} (i),
		\[
		\|\vp\|_\infty^{p-1}=\frac{p+1}{2}J_p^2=\Phi(\varepsilon_p),
		\qquad
		\Phi(\varepsilon):=\frac{\varepsilon}{2}\,B\left(\varepsilon,\frac12\right)^2.
		\]
		The function $\Phi$ is smooth and positive in $(0,+\infty)$, so that
		$p\mapsto\log\Phi(\varepsilon_p)$ is smooth in $(-1,+\infty)$. Since
		$B(\frac12,\frac12)=\Gamma(\frac12)^2=\pi$, we have $\Phi(\frac12)=\frac{\pi^2}{4}$.
		Moreover, using $\partial_a\log B(a,b)=(\log\Gamma)'(a)-(\log\Gamma)'(a+b)$,
		\[
		(\log\Phi)'(\varepsilon)=\frac1\varepsilon
		+2\left((\log\Gamma)'(\varepsilon)-(\log\Gamma)'\left(\varepsilon+\tfrac12\right)\right),
		\]
		and therefore, recalling that $(\log\Gamma)'(\frac12)-(\log\Gamma)'(1)=-2\log2$
		(see the proof of Proposition \ref{prop:new_as_behav_norm}),
		\[
		(\log\Phi)'\left(\tfrac12\right)=2-4\log2.
		\]
		Since $\varepsilon_1=\frac12$ and $\frac{d\varepsilon_p}{dp}\big|_{p=1}=-\frac14$,
		we get
		\[
		\frac{d}{dp}\log\Phi(\varepsilon_p)\Big|_{p=1}=-\frac14\,(2-4\log2)=\log2-\frac12.
		\]
		A Taylor expansion at $p=1$ then gives
		\[
		\log\|\vp\|_\infty^{p-1}=\log\frac{\pi^2}{4}+\left(\log2-\frac12\right)(p-1)
		+O\bigl((p-1)^2\bigr),
		\]
		and \eqref{esp_p_to_1} follows by exponentiating. Finally, comparing
		\eqref{esp_p_to_1} with Proposition \ref{Teo_p_to_1} (1), we obtain
		$\tilde c=\log2-\frac12$.
	\end{proof}
	
	\begin{osservazione}\label{rem:rescaling}
		Since, as can be seen by simple computations, for any $T\in\R$, $L>0$, the positive solution to $-v''=v^p$, with Dirichlet boundary conditions in $(-L+T,L+T)$, is
		\begin{equation}\label{rescaling_solution}
			\hat v_{p,L,T}(t):=\left(\frac1L\right)^{\frac{2}{p-1}}v_p\left(\frac{t-T}{L}\right),\quad t\in [-L+T,L+T],
		\end{equation}
		thus the previous results allow to deduce the asymptotic behaviour of $\hat v_{p,L,T}$ as $p\to+\infty$ or $p\to1$. 
	\end{osservazione}
	
	\section{Results on nodal solutions in $I$}\label{section:nodal}
	In this section we study the solution $\vpn$ of \eqref{problema_in_I}
	having $n$ nodal regions. In Subsection \ref{subsection:nodal_construction}
	we show that it is obtained by gluing $n$ alternating, translated and
	rescaled copies of the positive solution $\vp$ (see Proposition
	\ref{prop:existence_uniqueness}). This explicit structure allows us to
	transfer to $\vpn$ the results of Section \ref{section:positive}. In
	particular, Theorem \ref{teo:globale_any_n} follows directly from
	Propositions \ref{prop:monotonia_L_infty} and \ref{prop:ordine_potenziale},
	and its formula (1)
	is used in the subsequent asymptotic analysis. In Subsection
	\ref{Section:Nondeg_asymt_behav} we prove Theorems
	\ref{Teo_p_to_inf_rescaled} and \ref{Teo_p_to_1_rescaled}. The asymptotic
	behaviour of $\vpn$ follows by rescaling from Propositions
	\ref{Teo_p_to_inf}, \ref{prop:new_as_behav_norm}, \ref{Teo_p_to_1} and
	\ref{prop:explicit_c}. The asymptotics of the first eigenvalue
	$\tilde\alpha_1(p,n)$, instead, requires an additional argument.
	
	\subsection{Existence, uniqueness and properties of nodal solutions in $I$}\label{subsection:nodal_construction}
	
	We first characterize the nodal solution introduced above.
	
	\begin{proposizione} \label{prop:existence_uniqueness}
		Let $p >1$ and $n \geq 1$. Let $v_{p,n,i}:=\hat v_{p,\frac1n,s_{n,i}}$ in $I_{n,i}$ (where $\hat v_{p,\frac1n,s_{n,i}}$ is defined in \eqref{rescaling_solution}), namely
		\begin{equation}\label{positive_subsol}
			\vpin(t)=n^{\frac{2}{p-1}}\vp(n(t+1)+1-2i), \quad t \in I_{n,i},  \; i=1, \ldots, n. 
		\end{equation}
		Then the function 
		\begin{equation} \label{sum_solution}
			\vpn=\sum_{i=1}^{n}(-1)^{n+i}\barvpin,
		\end{equation}
		where $\barvpin$ denotes the trivial extension to the entire interval $I$ of $\vpin$, is the unique (up to sign) solution to \eqref{problema_in_I} having $n$ nodal regions. In particular, since $\bar v_{p,n,n}$ is taken with positive sign, $\vpn$ satisies
		$(v_p^{(n)})'(1)<0$.
	\end{proposizione}
	\begin{osservazione}
		Before proving the Proposition let us observe that if $n=1$
		\begin{equation}\label{vp1}
			v_p^{(1)}=v_{p,1,1}\equiv \vp,
		\end{equation}
		the positive solution to \eqref{positive_sol}.
	\end{osservazione}
	
	\begin{proof}
		In view of the above Remark, we focus on the case $n\geq2$.\\
		By construction, in each interval $I_{n,i}$ the function $\vpn$ coincides,
		up to sign, with $v_{p,n,i}$ and therefore satisfies the equation in
		\eqref{problema_in_I}. \\
		Moreover, at every nodal point $r_{n,i}$,
		$i=1,\ldots,n-1$, we have
		$v_{p,n,i}(r_{n,i})=v_{p,n,i+1}(r_{n,i})=0$ and
		$v_{p,n,i}'(r_{n,i})=-v_{p,n,i+1}'(r_{n,i})$.
		Hence, taking into account the alternating signs in \eqref{sum_solution},
		the function $\vpn$ is of class $C^1(I)$.
		
		Furthermore, since each $v_{p,n,i}$ satisfies
		$-v_{p,n,i}''=v_{p,n,i}^p$ in $I_{n,i}$ and vanishes at the endpoints
		of $I_{n,i}$, we have
		$v_{p,n,i}''(r_{n,i-1})=v_{p,n,i}''(r_{n,i})=0$.
		Thus the second derivatives also match at the nodal points, and
		$\vpn\in C^2(I)$. Therefore $\vpn$ solves \eqref{problema_in_I} in the
		whole interval $I$.
		
		The uniqueness up to sign follows from \cite[Theorem 6]{Beres}, where
		uniqueness is proved for a more general class of nonlinearities.
	\end{proof}

	\begin{osservazione}\label{rem:s_i}
		Since $v_{p,n,i}$ is a rescaling of the positive solution $v_p$ in $I$, it is immediate from what is known about $v_p$ that $v_{p,n,i}$ is symmetric with respect to $s_{n,i}$ and
		\[
		\vpin(s_{n,i})=\max_{I_{n,i}}\vpin.
		\]
		Furthermore, since $\vpn$ is a gluing of alternating copies of $v_{p,n,i}$, we also have
		$\|\vpn\|_{\infty}=\|v_{p,n,i}\|_{\infty}$
		for every $i\in\{1,\ldots,n\}$, and
		$v'_{p,n,i}(t)=0$
		if and only if $t=s_{n,i}$.
	\end{osservazione}
	
	We are now in position to prove the global properties stated in Theorem \ref{teo:globale_any_n}.
	
	\begin{proof}[Proof of Theorem \ref{teo:globale_any_n}]
		By Proposition \ref{prop:existence_uniqueness}, in each nodal interval
		$I_{n,i}$ the function $\vpn$ coincides, up to sign, with $\vpin$.
		Hence, by \eqref{positive_subsol},
		\begin{equation}\label{rescaling_nodal}
			|\vpn(t)|=n^{\frac{2}{p-1}}
			v_p\bigl(n(t-s_{n,i})\bigr),
			\qquad t\in I_{n,i}.
		\end{equation}
		
		Taking the maximum in \eqref{rescaling_nodal}, we obtain
		$\|\vpn\|_\infty=n^{\frac{2}{p-1}}\|v_p\|_\infty$. Therefore, by
		Proposition \ref{prop:monotonia_L_infty} (i),
		\[
		\|\vpn\|_\infty^{p-1}
		=
		n^2\|v_p\|_\infty^{p-1}
		=
		\frac{n^2(p+1)}{2}J_p^2,
		\]
		which proves (1).
		
		Let now $p_2>p_1>1$. If $t\in I_{n,i}$, then
		$n(t-s_{n,i})\in I$ and, by Proposition
		\ref{prop:ordine_potenziale},
		\[
		v_{p_2}\bigl(n(t-s_{n,i})\bigr)
		<
		v_{p_1}\bigl(n(t-s_{n,i})\bigr).
		\]
		Moreover,
		$n^{\frac{2}{p_2-1}}\leq n^{\frac{2}{p_1-1}}$. Hence, using
		\eqref{rescaling_nodal},
		\[
		\begin{aligned}
			|v_{p_2}^{(n)}(t)|
			&=n^{\frac{2}{p_2-1}}
			v_{p_2}\bigl(n(t-s_{n,i})\bigr)\\
			&<n^{\frac{2}{p_1-1}}
			v_{p_1}\bigl(n(t-s_{n,i})\bigr)
			=|v_{p_1}^{(n)}(t)|.
		\end{aligned}
		\]
		Thus (2) follows. In particular,
		$p\mapsto\|\vpn\|_\infty$ is strictly decreasing in $(1,+\infty)$.
		
		Finally, raising \eqref{rescaling_nodal} to the power $p-1$, we get
		\[
		p|\vpn(t)|^{p-1}
		=
		n^2p\,v_p\bigl(n(t-s_{n,i})\bigr)^{p-1}.
		\]
		Let $\tau=\tau(p_1,p_2)$ be given by Proposition
		\ref{prop:ordine_potenziale}. Since
		$|n(t-s_{n,i})|=n|t-s_{n,i}|$, the same proposition yields
		\[
		p_2|v_{p_2}^{(n)}(t)|^{p_2-1}
		>
		p_1|v_{p_1}^{(n)}(t)|^{p_1-1}
		\qquad\text{if }\quad
		|t-s_{n,i}|<\frac{\tau}{n},
		\]
		\[
		p_2|v_{p_2}^{(n)}(t)|^{p_2-1}
		<
		p_1|v_{p_1}^{(n)}(t)|^{p_1-1}
		\qquad\text{if }\quad
		\frac{\tau}{n}<|t-s_{n,i}|<\frac1n
		\]
		and the two potentials intersect exactly at
		$s_{n,i}-\frac{\tau}{n}$ and $s_{n,i}+\frac{\tau}{n}$ in each nodal
		interval. This proves (3).
	\end{proof}
	
	\subsection{Asymptotic behaviour of solutions}\label{Section:Nondeg_asymt_behav}
	Relying on Propositions \ref{Teo_p_to_inf} and \ref{Teo_p_to_1},
	we derive the asymptotic behaviour of the nodal solutions $\vpn$, defined in \eqref{sum_solution}, and of the first eigenvalue of \eqref{L_vpn} with Dirichlet boundary conditions in $I$.\\
	Before doing so, we make an observation that will be useful in what follows.
	\begin{osservazione}
		Being $\tilde\alpha_1(p,n)$ the first eigenvalue of \eqref{L_vpn} with Dirichlet boundary conditions in $I$, we can immediately observe that we have the following inequality    \begin{equation}\label{bound_alpha_fbelow}
			\tilde\alpha_1(p,n)=\inf_{q\in H^1_0(I)}\frac{\int_{-1}^{1}|q'(t)|^2\,dt-p\int_{-1}^{1}|\vpn(t)|^{p-1}q^2(t)\,dt}{\int_{-1}^{1}q^2(t)\,dt} \geq  -p\|\vpn\|_\infty^{p-1}.
		\end{equation}
	\end{osservazione} 
	Moreover we present a lemma that follows from the last part of the proof of Proposition 4.5 in \cite{GGPS}.
	\begin{lemma}\label{lemma:tildeK}
		Let $\eta>0$, $W$ the function defined in \eqref{W} and $\tilde K$ a solution to
		\begin{equation}\label{eq:tildeK}
			\begin{cases}
				-\tilde{K}''=e^W\tilde{K}-{\eta} \tilde{K}\qquad\text{in $\R$}\\
				\tilde{K}>0\\
				\|\tilde K\|_\infty=1,
			\end{cases}
		\end{equation}
		then $\eta=\frac12$.
	\end{lemma}

	We now prove the asymptotic results stated in Theorems \ref{Teo_p_to_inf_rescaled} and \ref{Teo_p_to_1_rescaled}.
	
	\begin{proof}[Proof of Theorem \ref{Teo_p_to_inf_rescaled}]
		We first prove \eqref{item:norm_p-1_inf} and \eqref{item:norm_inf}
		for every $n\geq1$. By Theorem \ref{teo:globale_any_n} (1) and
		Proposition \ref{prop:monotonia_L_infty} (i),
		$\|\vpn\|_\infty^{p-1}=n^2\|\vp\|_\infty^{p-1}$, hence
		\eqref{item:norm_p-1_inf} follows from \eqref{migliore_esp_p_to_infty}.
		Moreover, by Theorem \ref{teo:globale_any_n} (1),
		\[
		\log\|\vpn\|_\infty=\frac1{p-1}\left(\log\frac{n^2(p+1)}{2}+2\log J_p\right),
		\]
		and, by \eqref{logJp_expansion}, the bracket equals
		$\log\frac{n^2p}{2}+\frac{1+4\log2}{p}+O\left(\frac1{p^2}\right)$, namely
		\eqref{bracket_expansion} with $\log\frac p2$ replaced by $\log\frac{n^2p}{2}$.
		Since $\log\frac{n^2p}{2}=O(\log p)$, the computation following
		\eqref{bracket_expansion} in the proof of Proposition
		\ref{prop:new_as_behav_norm} applies verbatim and yields
		\eqref{item:norm_inf}.
		
		If $n=1$, in view of \eqref{vp1}, items (3), (4) and (5) are contained
		in Proposition \ref{Teo_p_to_inf}. If $n\geq2$, (3) and (4) follow
		immediately using \eqref{positive_subsol} and \eqref{sum_solution},
		rescaling separately in each subinterval $I_{n,i}$ and applying
		Proposition \ref{Teo_p_to_inf}.\\
		
		In order to show (5), for $n\geq2$, we consider a first
		Dirichlet eigenfunction $w_p>0$ of \eqref{L_vp}, associated with the
		first eigenvalue $\tilde\alpha_1(p)$, namely
		\[
		-w_p''-pv_p^{p-1}w_p=\tilde\alpha_1(p)w_p
		\qquad\text{in }I,
		\]
		with $w_p(\pm1)=0$. Fix $i\in\{1,\ldots,n\}$ and define
		\[
		\phi_{p,n}(t)
		:=
		\begin{cases}
			w_p\bigl(n(t-s_{n,i})\bigr),
			& t\in I_{n,i},\\
			0,
			& t\in I\setminus I_{n,i}.
		\end{cases}
		\]
		Since $w_p(\pm1)=0$, we have
		$\phi_{p,n}\in H_0^1(I)$. Moreover, in $I_{n,i}$,
		\[
		|v_p^{(n)}(t)|^{p-1}
		=
		n^2v_p\bigl(n(t-s_{n,i})\bigr)^{p-1}.
		\]
		Therefore, by the variational characterization of
		$\tilde\alpha_1(p,n)$ and via the change of variables
		$x=n(t-s_{n,i})$, we obtain
		\begin{eqnarray}
			\tilde\alpha_1(p,n)&\leq&
			\frac{\int_I|\phi_{p,n}'|^2\,dt
				-p\int_I|v_p^{(n)}|^{p-1}\phi_{p,n}^2\,dt
			}{\int_I \phi_{p,n}^2\,dt}\nonumber\\
			&=&n^2\frac{\int_I|w_p'|^2\,dx-p\int_Iv_p^{p-1}w_p^2\,dx
			}{\int_I w_p^2\,dx}
			=n^2\tilde\alpha_1(p).
			\label{inequality_alpha}
		\end{eqnarray}
		
		Combining \eqref{as_behav_norm} and \eqref{item:asbehavalpha1p} of
		Proposition \ref{Teo_p_to_inf} with the definition \eqref{mu_p} of $\mu_p$,
		we have
		\[
		\mu_p^2\tilde\alpha_1(p)\to-\frac12.
		\]
		Moreover,
		\[
		\frac{\mu_p^2}{n^2}
		=
		\frac{1}{n^2p\|\vp\|_\infty^{p-1}}
		=
		\frac{1}{p\|\vpin\|_\infty^{p-1}}
		=
		\mupn^2.
		\]
		Hence, by \eqref{inequality_alpha},
		\[
		\mupn^2\tilde\alpha_1(p,n)
		=
		\frac{\mu_p^2}{n^2}\tilde\alpha_1(p,n)
		\leq
		\mu_p^2\tilde\alpha_1(p)
		\to-\frac12.
		\]
		On the other hand, by \eqref{bound_alpha_fbelow},
		\[
		\mupn^2\tilde\alpha_1(p,n)\geq-1.
		\]
		Thus, up to a subsequence,
		\[
		\mupn^2\tilde\alpha_1(p,n)
		\to
		\tilde\beta_1
		\in\left[-1,-\frac12\right].
		\]
		
		Finally, denoting by $\wpn$ the positive eigenfunction associated to
		$\tilde\alpha_1(p,n)$ satisfying $\|\wpn\|_\infty=1$, it solves
		\[
		-\wpn''
		=
		p|\vpn|^{p-1}\wpn+\tilde\alpha_1(p,n)\wpn
		\qquad\text{in }I.
		\]
		Let $\xi_p$ be such that
		\[
		\wpn(\xi_p)=\|\wpn\|_\infty=1.
		\]
		Up to a subsequence, we can assume that
		$\xi_p\in\overline{I}_{n,i}$ for some fixed
		$i\in\{1,\ldots,n\}$.
		
		We first observe that, for $p$ sufficiently large, $\xi_p\in I_{n,i}$.
		Indeed, $\tilde\alpha_1(p,n)<0$ for $p$ sufficiently large.
		Clearly $\xi_p\neq\pm1$, since $\wpn(\pm1)=0$. If $\xi_p$ were an
		interior nodal point of $\vpn$, then $\vpn(\xi_p)=0$ and, by the
		equation satisfied by $\wpn$,
		\[
		-\wpn''(\xi_p)
		=
		\tilde\alpha_1(p,n)\wpn(\xi_p)
		=
		\tilde\alpha_1(p,n)<0.
		\]
		On the other hand, since $\xi_p$ is a maximum point of $\wpn$,
		we have $\wpn''(\xi_p)\leq0$, and hence a contradiction.
		
		We define
		\[
		\wpntildei(y)
		:=
		\wpn(\mupn y+s_{n,i}),
		\qquad
		y\in
		\left[
		\frac{-1-s_{n,i}}{\mupn},
		\frac{1-s_{n,i}}{\mupn}
		\right].
		\]
		Since $\mupn=\mu_p/n$, the interval $\left[-\frac1{\mu_p},\frac1{\mu_p}\right]$
		corresponds precisely to the nodal interval
		$\overline I_{n,i}$. Moreover, by \eqref{positive_subsol},
		$\vpntildei=\tilde v_p$ in $
		\left[-\frac1{\mu_p},\frac1{\mu_p}\right]$,
		where $\tilde v_p$ is defined in \eqref{utildep}. Therefore
		\begin{equation}\label{eq:wtilde1p}
			-\wpntildei''
			=
			\left(1+\frac{\vpntildei}{p}\right)^{p-1}\wpntildei
			+\mupn^2\tilde\alpha_1(p,n)\wpntildei
			\qquad
			\text{in }
			\left(-\frac1{\mu_p},\frac1{\mu_p}\right).
		\end{equation}
		Observing that $\vpntildei(y)\in[-p,0]$ and using
		$\|\wpntildei\|_\infty\leq1$, by \eqref{eq:wtilde1p} we have,
		for $p$ sufficiently large,
		\[
		|\wpntildei''(y)|
		\leq
		1+\mupn^2|\tilde\alpha_1(p,n)|
		\leq2.
		\]
		We now show that $\wpntildei'$ is uniformly bounded on compact
		subsets of $\mathbb R$. Fix $R>0$. Since $\mu_p\to0$, for $p$
		sufficiently large, $
		[-R-1,R+1]
		\subset
		\left(-\frac1{\mu_p},\frac1{\mu_p}\right)$.
		For every $y\in[-R,R]$, by the mean value theorem there exists
		$z_p\in(y,y+1)$ such that
		\[
		|\wpntildei'(z_p)|
		=
		|\wpntildei(y+1)-\wpntildei(y)|
		\leq2.
		\]
		Hence, using the bound on $\wpntildei''$,
		\[
		|\wpntildei'(y)|
		\leq
		|\wpntildei'(z_p)|
		+
		\int_y^{z_p}|\wpntildei''(s)|\,ds
		\leq4.
		\]
		Thus, for every $R>0$,
		$\|\wpntildei'\|_{L^\infty(-R,R)}
		\leq C_R$, therefore, by the Ascoli--Arzelà theorem, up to a subsequence,
		\[
		\wpntildei\to\psi_{n,i}
		\qquad\text{in }C^1_{\rm loc}(\mathbb R).
		\]
		By \eqref{mu_p_piu} of Proposition \ref{Teo_p_to_inf}, $\psi_{n,i}$ satisfies
		\[
		\begin{cases}
			-\psi_{n,i}''
			=
			e^W\psi_{n,i}
			+\tilde\beta_1\psi_{n,i}
			&\text{in }\mathbb R,\\
			\psi_{n,i}\geq0.
		\end{cases}
		\]
		Set
		\[
		\tilde\xi_p:=\frac{\xi_p-s_{n,i}}{\mupn}.
		\]
		Since $\xi_p\in I_{n,i}$, $\tilde\xi_p
		\in
		\left(-\frac1{\mu_p},\frac1{\mu_p}\right)$,
		and, since $\xi_p$ is a maximum point of $\wpn$,
		\[
		\wpntildei(\tilde\xi_p)=1,
		\qquad
		\wpntildei'(\tilde\xi_p)=0,
		\qquad
		\wpntildei''(\tilde\xi_p)\leq0.
		\]
		
		We claim that $(\tilde\xi_p)_p$ is bounded. Suppose by contradiction
		that, along a subsequence,  $|\tilde\xi_p|\to+\infty$.
		Then, by \eqref{eq:wtilde1p}, \eqref{W} and (5) of Proposition $\ref{Teo_p_to_inf}$, for $p$
		sufficiently large,
		\[
		0\leq
		-\wpntildei''(\tilde\xi_p)
		\leq
		e^{\frac{p-1}{p}
			(C_1W(\tilde\xi_p)+C_2)}
		+\tilde\beta_1+o_p(1)
		=
		\tilde\beta_1+o_p(1)<0,
		\]
		which is impossible. Hence $
		|\tilde\xi_p|\leq C$,
		and, up to a subsequence,  $
		\tilde\xi_p\to\tilde\xi_\infty\in\mathbb R$.
		By the $C^1_{loc}$ convergence of $\wpntildei$ to $\psi_{n,i}$,
		\[
		1
		=
		\wpntildei(\tilde\xi_p)
		\longrightarrow
		\psi_{n,i}(\tilde\xi_\infty).
		\]
		Therefore $\psi_{n,i}(\tilde\xi_\infty)
		=1=\|\psi_{n,i}\|_\infty$. From the uniqueness for the corresponding Cauchy problem we can also deduce that $\psi_{n,i}>0$.
		Summing up, $\psi_{n,i}$ solves \eqref{eq:tildeK}, so applying Lemma \ref{lemma:tildeK} we have that $\tilde\beta_1=-\frac12$.\\
		The thesis follows from \eqref{mupn} and \eqref{item:norm_p-1_inf}.
	\end{proof}

	\begin{proof}[Proof of Theorem \ref{Teo_p_to_1_rescaled}]
		By Theorem \ref{teo:globale_any_n} (1) and Proposition
		\ref{prop:monotonia_L_infty} (i), $\|\vpn\|_\infty^{p-1}=n^2\|\vp\|_\infty^{p-1}$,
		hence (1) follows from Proposition \ref{prop:explicit_c}.
		If $n=1$, in view of \eqref{vp1}, (2), (3) and (4) are contained in
		Proposition \ref{Teo_p_to_1}. If $n\geq2$, (3) can be immediately obtained
		using \eqref{sum_solution} and \eqref{positive_subsol}, rescaling separately
		in each subinterval $I_{n,i}$ and applying Proposition \ref{Teo_p_to_1}.\\
		Let us now show (2) for $n\geq2$. First of all by \eqref{sum_solution} and \eqref{positive_subsol}
		\begin{equation}\label{A}
			\frac{\vpn}{\|\vpn\|_\infty}
			=\sum_{i=1}^n(-1)^{n+i}
			\frac{\barvpin}{\|\vpin\|_\infty}.
		\end{equation}
		Next \eqref{positive_subsol} and (2) of Theorem $\ref{Teo_p_to_1}$ imply that \text{as $p\to1^+$}
		\begin{equation}\label{B}
			\frac\vpin{\|\vpin\|_\infty}\to\varphi_{1,n,i}\qquad\text{in $C^1(\bar I_{n,i})$},
		\end{equation}
		where $\varphi_{1,n,i}(t):=\varphi_1(n(t+1)+1-2i)$.\\In conclusion,
		\[
		\sum_{i=1}^n(-1)^{n+i}\overline{\varphi}_{1,n,i}
		\]
		is an eigenfunction of $-\frac{d^2}{dt^2}$ in $I$ corresponding to
		$\nu_n$, with $L^\infty$ norm equal to $1$. Moreover,
		\[
		\left(\sum_{i=1}^n(-1)^{n+i}\overline{\varphi}_{1,n,i}\right)'(1)=\varphi_{1,n,n}'(1)=n\varphi_1'(1)<0.
		\]
		Hence
		\[
		\varphi_n=\sum_{i=1}^n(-1)^{n+i}\overline{\varphi}_{1,n,i},
		\]
		and (2) follows from \eqref{A} and \eqref{B}.\\
		Lastly let us prove (4) for $n\geq2$. Combining \eqref{bound_alpha_fbelow} and (1) of Theorem $\ref{Teo_p_to_1_rescaled}$ we have
		\[
		\tilde\alpha_1(p,n)\geq -n^2\frac{\pi^2}{4}+o_p(1),
		\] 
		hence, being $\tilde\alpha_1(p,n)<0$ (as it can be seen testing the Rayleigh quotient defining $\tilde\alpha_1(p,n)$ on $\vpn$), $\tilde\alpha_1(p,n)$ is bounded and there exists $\tilde\alpha\in[-n^2\frac{\pi^2}{4},0]$ such that, up to a subsequence, $\tilde\alpha_1(p,n)\to\tilde\alpha$ as $p\to 1^+$.\\
		Let $\wpn$ be the positive eigenfunction of \eqref{L_vpn} associated to $\tilde\alpha_1(p,n)$ and satisfying $\|\wpn\|_\infty=1$. Therefore, by (1),
		\[
		|\wpn''|\leq (p|\vpn|^{p-1}+|\tilde\alpha_1(p,n)|)|\wpn|\leq C\quad\text{and in turn}\quad |\wpn'|\leq 2C,
		\]
		thus up to a subsequence $\wpn\to\bar w$ in $C^1(\bar I)$, with $\bar w\geq0$ in $I$ and $\bar w(\pm1)=0$. Moreover, by (3),
		$p|\vpn|^{p-1}\to n^2\frac{\pi^2}{4}$ a.e. in $I$, and the sequence
		$p|\vpn|^{p-1}$ is uniformly bounded in $L^\infty(I)$. Hence, passing
		to the limit in the weak formulation, we obtain
		\[
		-\bar w''=
		\left(n^2\frac{\pi^2}{4}+\tilde\alpha\right)\bar w
		\quad\text{in }I
		\]
		in the weak sense. Since $\bar w\in C^1(\bar I)$, the right-hand side is
		continuous, and therefore $\bar w\in C^2(I)$ and the equation holds
		classically.
		Since $\wpn\to\bar w$ uniformly in $\bar I$ and $\|\wpn\|_\infty=1$, we get
		$\|\bar w\|_\infty=\lim_{p\to1^+}\|\wpn\|_\infty=1$, so that $\bar w\not\equiv0$.
		Then $\bar w>0$ and hence $\bar w$ is an eigenfunction associated to the first eigenvalue $\nu_1=\frac{\pi^2}{4}$ of $-\frac{d^2}{dt^2}$ in $I$ with Dirichlet boundary conditions, thus $n^2\frac{\pi^2}{4}+\tilde\alpha=\frac{\pi^2}{4}$, and so
		\[
		\tilde\alpha=\left(1-n^2\right)\frac{\pi^2}{4}.
		\]
		Since every converging subsequence of $\tilde\alpha_1(p,n)$ has the same
		limit, the whole sequence converges, and (4) follows.
	\end{proof}
	
	\section{Stability and instability results}\label{section:stability}
	
	In this section we recall the variational framework introduced in
	\cite{AfonsoPacella} and we prove Theorems \ref{Teo_stability_p_infty}
	and \ref{Teo_stability_p_1}. In Subsection \ref{subsection:energy} we
	define the local shape functional associated with a nondegenerate
	solution and the notions of stable and unstable energy-stationary pairs
	under the volume constraint. In Subsection \ref{subsection:one-dim sol}
	we focus on flat cylinders: we recall a sufficient condition for the
	nondegeneracy of the one-dimensional solutions and the characterization
	of their stability given in \cite[Theorem 4.11]{AfonsoPacella} (see
	Theorem \ref{teorema}), which reduces the problem to the sign of the
	derivative at the endpoint of the solution of an auxiliary linear
	ordinary differential equation, and we observe that, after rescaling,
	it is enough to consider cylinders of height one. Finally, in
	Subsections \ref{subsection:proof_infty} and \ref{subsection:proof_1}
	we combine this characterization with the asymptotic results of
	Section \ref{section:nodal}.

	\subsection{Energy functional and energy stationary pairs}\label{subsection:energy}
	Let $\omega\subset\R^{N-1}$ be a smooth bounded domain and let $\mathcal C_\omega$ be the half-cylinder spanned by $\omega$, namely $\mathcal C_\omega:=\omega\times(0,+\infty).
	$
	We denote by $x=(x',x_N)$ the points in $\bar{\mathcal C}_\omega$, where $x'=(x_1,\ldots,x_{N-1})\in\bar\omega$ and $x_N\geq0$.\\
	In $\mathcal C_\omega$ we consider domains whose relative boundaries are cartesian graphs of functions in $C^2(\bar\omega)$. More precisely, for $\varphi\in C^2(\bar\omega)$, we set
	\begin{align*}
		&\Omega_\varphi:=\{(x',x_N)\in\mathcal C_\omega\,:\,x_N<e^{\varphi(x')}\},\\
		&\Gamma_\varphi:=\{(x',x_N)\in\mathcal C_\omega\,:\,x_N=e^{\varphi(x')}\},\\
		&\Gamma_{1,\varphi}:=(\partial \Omega_\varphi\setminus \bar\Gamma_\varphi).
	\end{align*}
	
	We will consider variations of $\Omega_\varphi$ in the class of cartesian graphs of the type $\Omega_{\varphi+t\psi}$ for $\psi \in C^2(\bar\omega)$, which amounts to consider a one parameter family of diffeomorphisms $\xi:(-\eta,\eta)\times\bar{\mathcal C}_\omega\to\bar{\mathcal C}_\omega$ of the type
	\[
	\xi(t,x)=(x',e^{t\psi (x')}x_N),
	\]
	which is generated by the vector field $V(x)=(0',\psi (x')x_N)$, where $0'=(0,\ldots,0)\in\R^{N-1}$.\\
	Let $\bar\varphi\in C^2(\bar\omega)$ and let $u_{\Omega_{\bar\varphi}}\in W^{1,\infty}(\Omega_{\bar\varphi})\cap W^{2,2}(\Omega_{\bar\varphi})$ be a nondegenerate solution to \eqref{problem_Omega}, with $\Omega=\Omega_{\bar\varphi}$. In \cite[Proposition 2.1 and Remark 2.8]{AfonsoPacella} it has been shown that under such deformations of $\Omega_{\bar\varphi}$ the non-degeneracy of $u_{\Omega_{\bar\varphi}}$ induces a local uniqueness result for solutions of \eqref{problem_Omega} in the deformed domains. Namely, given $\psi \in C^2(\bar\omega)$, there exists $\delta>0$ such that for any $t\in(-\delta,\delta)$ problem \eqref{problem_Omega} with $\Omega=\Omega_{\bar\varphi+t\psi}$ admits a unique solution $u_{\Omega_{\bar\varphi+t\psi}}$ in a neighborhood of $u_{\Omega_{\bar\varphi}}\circ\xi(t,\cdot)^{-1}$. 
	Thus, for any $\psi\in C^2(\bar\omega)$, the energy functional 
	\[
	T(\Omega_{\bar\varphi+t\psi})=J(u_{\Omega_{\bar\varphi+t\psi}})=\frac12\int_{\Omega_{\bar\varphi+t\psi}}|\nabla u_{\Omega_{\bar\varphi+t\psi}}(x)|^2\,dx-\frac{1}{p+1}\int_{\Omega_{\bar\varphi+t\psi}} |u_{\Omega_{\bar\varphi+t\psi}}|^{p+1}(x)\,dx
	\]
	is well defined for $t$ sufficiently small.\\
	Ultimately the energy functional $T$ is locally well defined as a functional on $C^2(\bar\omega)$ near $\bar\varphi$. With an abuse of notation, for any $\psi\in C^2(\bar\omega)$ and for $t$ sufficiently small we set
	\[
	T(\bar\varphi+t\psi):=T(\Omega_{\bar\varphi+t\psi}).
	\]
	Moreover in \cite[Lemma 4.1]{AfonsoPacella} the first derivative of $T$ at $\bar\varphi$, i.e. for $t=0$, with respect to variations $\psi \in C^2(\bar\omega)$ has been computed and it takes the following form
	\begin{equation}\label{T'}
		T'(\bar\varphi)[\psi]=-\frac12\int_\omega \left(\frac{\partial u_{\Omega_{\bar\varphi}}}{\partial \nu}(x',e^{\bar\varphi(x')})\right)^2 \psi(x') e^{\bar\varphi(x')}dx'.
	\end{equation} 
	We point out that all the results in \cite[Section 4]{AfonsoPacella} hold also for non-degenerate sign-changing solutions of \eqref{problem_Omega} with $\Omega=\Omega_\varphi$.\\
	
	We will be interested in understanding how the energy of a solution $u_{\Omega_{\varphi}}$ behaves with respect to volume-preserving variations of $\Omega_{\varphi}$.\\
	First of all we notice that the volume functional 
	\[
	\mathcal V(\varphi):=\int_\omega e^{\varphi(x')}dx'
	\]
	is of class $C^2$ and $\mathcal V'(\varphi)[\psi]=\int_\omega e^{\varphi(x')} \psi(x') \,dx'$ for any $\psi\in C^2(\bar\omega)$.\\
	Consequently we define the manifold 
	\[
	M:=\{\varphi\in C^2(\bar\omega)\,:\,\int_{\omega}e^{\varphi (x')}dx'=\int_{\omega}e^{\bar\varphi(x')}dx'\}
	\]
	and we consider the restricted functional
	\[
	I(\varphi):=T_{|M}(\varphi),\quad{\varphi\in M}.
	\]
	Obviously $\bar\varphi\in M$ and the tangent space at $\bar\varphi$ is \[
	T_{\bar\varphi}M=\{\psi\in C^2(\bar\omega)\,:\,\int_\omega e^{\bar\varphi(x')}\psi(x')\,dx'=0\}.
	\]
	\begin{definizione}\label{def:critical_energy_pair}
		We say that $(\Omega_{\bar\varphi},u_{\Omega_{\bar\varphi}})$ is an energy-stationary pair under volume constraint if $\bar\varphi$ is a critical point of $I$, namely if $I'(\bar\varphi)[\psi]=0$ for any $\psi\in T_{\bar\varphi}M$, or equivalently if there exists $\lambda\in\R$ such that $T'(\bar\varphi)=\lambda \mathcal V'(\bar\varphi)$.
	\end{definizione}
	
	\begin{definizione}\label{def:critical_energy_pair_stable}
		An energy-stationary pair $(\Omega_{\bar\varphi},u_{\Omega_{\bar\varphi}})$ is called stable if 
		\[
		I''(\bar\varphi)[\psi,\psi]>0\qquad\forall\: \psi\in T_{\bar \varphi}M\setminus\{0\}.
		\]
		where $I''$ denotes the second derivative of the restricted functional $I$
	\end{definizione}
	
	\begin{definizione}\label{def:critical_energy_pair_unstable}
		An energy-stationary pair $(\Omega_{\bar\varphi},u_{\Omega_{\bar\varphi}})$ is called unstable if there exists $\psi\in T_{\bar \varphi}M$ such that
		\[
		I''(\bar\varphi)[\psi,\psi]<0.
		\]
	\end{definizione}
	Clearly if $(\Omega_{\bar\varphi},u_{\Omega_{\bar\varphi}})$ is an unstable energy pair, then $\bar\varphi$ is not a local minimizer for $T$ under volume constraint.
	\subsection{One-dimensional solutions}\label{subsection:one-dim sol}
	Let $L>0$, let $\varphi_L(x'):=\log L$, $x'\in\bar\omega$, and let us consider 
	\[
	\Omega_{\varphi_L}=\{(x',x_N)\in \mathcal C_\omega\,:\,x_N< L\}.
	\]
	When $\Omega=\Omega_{\varphi_L}$, any weak solution of \eqref{problem_Omega} is also a classical solution up to the boundary, i.e. it belongs to $C^2(\overline{\Omega_{\varphi_L}})$. This follows by standard regularity theory by considering the boundary conditions and that $\partial\Omega_{\varphi_L}$ is made by the union of three $(N-1)$-dimensional manifolds (with boundary) intersecting orthogonally (see also \cite[Proposition 6.1]{PacellaTralli}).\\
	In this subsection we will introduce and discuss properties of one-dimensional solutions to \eqref{problem_Omega}, when $\Omega=\Omega_{\varphi_L}$.\\
	
	Let $\upLm$ be the solution of \eqref{ODE_in_(0,L)} having $m$ nodal regions and verifying $(\upLm)'(L)<0$. We denote by $\vpLdmu$ the even extension of $\upLm$ to $(-L,L)$, which solves the corresponding Dirichlet problem in $(-L,L)$ and has $2m-1$ nodal regions. Moreover, with an abuse of notation, we denote by $\upLm(x',x_N):=\upLm(x_N)$ the one-dimensional solution of \eqref{problem_Omega} in $\Omega_{\varphi_L}$.
	
	
	Sufficient conditions for nondegeneracy of $\upLm$ (as a solution to \eqref{problem_Omega}) have been obtained by decomposing the spectrum of the linearized operator $-\Delta-p|\upLm|^{p-1}$ in the space $H^1_0(\Omega_{\varphi_L}\cup \Gamma_{1,\varphi_L})$.
	
	\begin{proposizione}\cite[Corollary 4.7]{AfonsoPacella}\label{prop:nondegeneracy}
		Let $p>1$, $m\geq1$, then $\upLm$ is a nondegenerate solution of \eqref{problem_Omega} in $\Omega_{\varphi_L}$ if both the following conditions are satisfied:
		\begin{itemize}
			\item[\emph{i.}] the eigenvalue problem 
			\begin{equation}\label{eigenvalue_problem_alpha_L}
				\begin{cases}
					-z'' - p|\upLm|^{p-1}z=\alpha z \quad \mbox{ in } (0,L)\\
					z'(0)=z(L)=0
				\end{cases}
			\end{equation}
			does not admit zero as an eigenvalue;
			\item[\emph{ii.}] $\lambda_1(\omega)>-{\alpha}_{1,L}(p,m)$, where ${\alpha}_{1,L}(p,m)$ is the first eigenvalue of \eqref{eigenvalue_problem_alpha_L}.
		\end{itemize}   
	\end{proposizione}
	\begin{osservazione}\label{rem:nondegeneracy}
		Being $\alpha_j(p,m)=L^2\alpha_{j,L}(p,m)$ for any $j$, \cite[Proposition 3.8]{Afonso1} implies that condition {i.} of Proposition \ref{prop:nondegeneracy} is always satisfied. As a consequence, if $\lambda_1(\omega)>-\alpha_{1,L}(p,m)$, the one-dimensional solution $\upLm$ to \eqref{problem_Omega} in $\Omega_{\varphi_L}$ is nondegenerate.   
	\end{osservazione}
	
	\begin{osservazione}\label{rem:tildealpha}
		We can observe that
		$\alpha_{1,L}(p,m)=\tilde\alpha_{1,L}(p,2m-1)$,
		where $\tilde\alpha_{1,L}(p,2m-1)$ is the first eigenvalue of the linearized operator \eqref{L_vpn} (for $n=2m-1$) in $(-L,L)$ with Dirichlet boundary conditions.
		Indeed, if
		$V=\{\varphi\in H^1(0,L):\varphi(L)=0\}$,
		then the even extension to $(-L,L)$ of any $\varphi\in V$ belongs to $H^1_0(-L,L)$, hence
		$\tilde\alpha_{1,L}(p,2m-1)\leq\alpha_{1,L}(p,m)$.
		Conversely, let $w>0$ be a first eigenfunction associated with
		$\tilde\alpha_{1,L}(p,2m-1)$. Since the potential
		$p|v_{p,L}^{(2m-1)}|^{p-1}$ is even and the first Dirichlet eigenvalue is simple, $w$ is even. Thus $w'(0)=0$, and the restriction of $w$ to $(0,L)$ is admissible for the variational characterization of $\alpha_{1,L}(p,m)$. Therefore
		$\alpha_{1,L}(p,m)\leq\tilde\alpha_{1,L}(p,2m-1)$,
		and the conclusion follows.
	\end{osservazione}

	Therefore, if we define
	\[
	M=\{\varphi\in C^2(\bar\omega)\,:\,\int_\omega e^{\varphi(x')}dx'=\int_\omega e^{\varphi_L(x')}dx'=L|\omega|\},
	\]
	then $(\Omega_{\varphi_L}, \upLm)$ is an energy-stationary pair under volume constraint, being, by \eqref{T'},
	\[
	T'(\varphi_L)[\psi]=-\frac{L}{2}((\upLm)'(L))^2\int_\omega \psi(x')\,dx'=0
	\]
	for any $\psi\in T_{\varphi_L}M=\{\psi\in C^2(\bar\omega)\,:\,\int_\omega \psi(x')\,dx'=0\}$.
	
	It is worth pointing out that assumption $\lambda_1(\omega)>-\alpha_{1,L}(p,m)$ is not only a sufficient condition to guarantee the non-degeneracy of the one-dimensional solution $\upLm$ of \eqref{problem_Omega} in $\Omega_{\varphi_L}$ (see Remark \ref{rem:nondegeneracy}), but also implies the positivity of the function $\hpLm$,
	which is a key property in the proof of Theorem \ref{teorema}. We state the positivity result, obtained in  \cite[Proposition 4.10]{AfonsoPacella} for $m=1$, in a suitable form and for $L=1$. 
	\begin{proposizione}\label{prop:hpL>0} 
		
		Let $\lambda>0$ and let assume that $\lambda>-\alpha_{1}(p,m)$, where $\alpha_{1}(p,m)$ is the first eigenvalue of \eqref{eigenvalue_problem_alpha_L} for $L=1$.\\
		Then 
		\begin{equation}\label{hptilde}
			\begin{cases}
				h''=(\lambda-p|\upm|^{p-1}) h \quad \mbox{ in } (0,1)\\
				h(1)=-(\upm)'(1)>0\\
				h'(0)=0.
			\end{cases}
		\end{equation}
		admits a unique solution, which is positive in $[0,1]$.
		\begin{proof}
			We first observe that problem \eqref{hptilde} admits a unique solution.
			Indeed, let $y$ be the solution of the Cauchy problem
			$y''=(\lambda-p|\upm|^{p-1})y$, with $y(0)=1$ and $y'(0)=0$.
			Then $y(1)\neq0$. Otherwise, $y$ would be a nontrivial solution of
			\[
			-y''-p|\upm|^{p-1}y=-\lambda y
			\quad\text{in }(0,1),
			\qquad
			y'(0)=0,\quad y(1)=0,
			\]
			and hence $-\lambda$ would be an eigenvalue of
			\eqref{eigenvalue_problem_alpha_L}, which is impossible since
			$-\lambda<\alpha_1(p,m)$.
			
			Therefore
			\[
			h(t):=-\frac{(\upm)'(1)}{y(1)}\,y(t)
			\]
			is the unique solution of \eqref{hptilde}.
			
			Let now $z_1>0$ be a first eigenfunction associated with
			$\alpha_1(p,m)$. Then
			$z_1'(0)=0$, $z_1(1)=0$, and
			\[
			z_1''+
			\bigl(p|\upm|^{p-1}+\alpha_1(p,m)\bigr)z_1=0.
			\]
			On the other hand,
			\[
			h''+
			\bigl(p|\upm|^{p-1}-\lambda\bigr)h=0,
			\]
			and, by assumption,
			\[
			p|\upm|^{p-1}-\lambda
			<
			p|\upm|^{p-1}+\alpha_1(p,m).
			\]
			
			Suppose by contradiction that $h$ has a zero in $(0,1)$, and let
			$\xi$ be its first positive zero. Since $h(0)\neq0$, we may assume
			that $h(0)>0$ (otherwise we can replace $h$ by $-h$ for the purpose
			of the comparison). Thus $h>0$ in $[0,\xi)$ and $h'(0)=0$.
			
			We extend both $h$ and $z_1$ evenly to $(-1,1)$. Then $-\xi$ and
			$\xi$ are two consecutive zeros of the even extension of $h$.
			Since the coefficient in the equation satisfied by $z_1$ is strictly
			larger than the one in the equation satisfied by $h$, the Sturm
			comparison theorem implies that the even extension of $z_1$ must have
			a zero in $(-\xi,\xi)$. This is impossible, since $z_1$ is a first
			eigenfunction and therefore its even extension is strictly positive
			in $(-1,1)$.
			
			Hence $h$ has no zero in $(0,1)$. Since $h(1)>0$, it follows that $
			h>0$ in $[0,1]$.
		\end{proof}
	\end{proposizione}

	Next we recall a result obtained in \cite{AfonsoPacella}, that will be crucial in order to prove Theorem \ref{Teo_stability_p_infty} and Theorem \ref{Teo_stability_p_1}.
	\begin{teorema}\label{teorema}
		Let $\omega\subset\R^{N-1}$ be a smooth bounded domain and let $\lambda_1(\omega)$ be the first nontrivial Neumann eigenvalue of $-\Delta_{\R^{N-1}}$ in $\omega$.\\
		Let $p>1$, $m\geq1$ and let $\upLm$ be the one-dimensional solution with $m$ nodal regions to \eqref{problem_Omega} in $\Omega_{\varphi_L}$, such that $\frac{\partial \upLm}{\partial\nu}<0$ on $\Gamma_{\Omega_{\varphi_L}}$.\\
		Let us assume that $\lambda_1(\omega)>-\alpha_{1,L}(p,m)$, where $\alpha_{1,L}(p,m)$ is the first eigenvalue of \eqref{eigenvalue_problem_alpha_L}, and let $\hpLm$ be the solution to 
		\begin{equation}\label{hpLkequation}
			\begin{cases}
				\hpLm''=(\lambda_1(\omega)-p|\upLm|^{p-1})\,\hpLm \quad \mbox{ in } (0,L)\\
				\hpLm(L)=-L\,(\upLm)'(L)\\
				\hpLm'(0)=0.
			\end{cases}
		\end{equation}
		Then:\\
		\emph{$i)$} $(\Omega_{\varphi_L}, \upLm)$ is an unstable energy stationary pair if and only if $\hpLm'(L)<0$;\\
		\emph{$ii)$} $(\Omega_{\varphi_L}, \upLm)$ is a stable energy stationary pair if and only if $\hpLm'(L)>0$. 
	\end{teorema}
	In fact, in \cite[Theorem 4.11]{AfonsoPacella} this result has been obtained in the case $L=1$, but it can be easily extended to the general case $L>0$.
	Furthermore, in the above-mentioned theorem, the positivity of the one-dimensional solution is assumed. However, this assumption is used only to ensure that the derivative of the one-dimensional solution at the right endpoint of the interval is negative. Therefore, the theorem also applies in our setting,  since we assume this property.
	
	More importantly, in the case $L=1$, the quantity $\rho$ appearing in the statement of \cite[Theorem 4.11]{AfonsoPacella} is, in our Lane-Emden setting, equal to $-u'(1)h'(1)$; see \cite[(4.33)]{AfonsoPacella}. In the case of the cylinder of height $L$, the corresponding quantity is
	\begin{equation}\label{rho}
		\rho=-L\,(\upLm)'(L)\hpLm'(L).
	\end{equation}
	Moreover, although the positivity and negativity of $\rho$ are stated there only as sufficient conditions for stability and instability, respectively, they actually hold as \virg{if and only if} conditions, as is clear from the proof of \cite[Theorem 4.11]{AfonsoPacella} (see also \cite[Theorem 1.5]{AfonsoPacella}). Indeed, the proof shows that
	\[
	I''(\varphi_L)[\psi_1,\psi_1]=
	\rho\quad\text{and}\quad I''(\varphi_L)[v,v]\geq
	\rho\,
	\qquad\text{for every }v\in T_{\varphi_L}M\;\text{s.t.}\;\|v\|_{L^2(\omega)}^2=1,
	\]
	where $\psi_1$ is an eigenfunction associated to the first non-trivial eigenvalue of the Laplacian in $\omega$ with Neumann boundary conditions, satisfying $\|\psi_1\|_{L^2(\omega)}^2=1$.
	In particular, if $\hpLm'(L)>0$, then $\rho>0$ and the second variation is uniformly positive, in the $L^2(\omega)$ sense, on $T_{\varphi_L}M$. Consequently, $(\Omega_{\varphi_L},\upLm)$ locally minimizes the energy along every nontrivial sufficiently small volume-preserving one-parameter perturbation. On the other hand, if $\hpLm'(L)<0$, then $\rho<0$ and $\psi_1\in T_{\varphi_L}M$ is a direction along which the second variation is negative.
	
	\begin{osservazione}\label{rem:casonegativo}
		The same stability/instability properties hold for $-\upLm$. Indeed, if $\bar h$ denotes the function associated with $-\upLm$, then
		\[
		\begin{cases}
			\bar h''=(\lambda_1(\omega)-p|\upLm|^{p-1})\bar h & \mbox{in }(0,L),\\
			\bar h(L)=L(\upLm)'(L),\\
			\bar h'(0)=0,
		\end{cases}
		\]
		hence, by uniqueness, $\bar h=-\hpLm$. Therefore $-L\,(-\upLm)'(L)\bar h'(L) = -L\,(\upLm)'(L)\hpLm'(L)$,
		so the quantity $\rho$, defined in \eqref{rho}, remains unchanged when $\upLm$ is replaced by $-\upLm$. Consequently, the pairs $(\Omega_{\varphi_L},\pm\upLm)$ have the same stability/instability properties: they are stable if $\hpLm'(L)>0$ and unstable if $\hpLm'(L)<0$.
	\end{osservazione}

	\begin{osservazione}\label{rem:equivalence}
		It is convenient for our purposes to define in $(0,1)$ a rescaling $\hpm$ of $\hpLm$, namely 
		\[
		\hpm(t):=L^{\frac{2}{p-1}}\hpLm(Lt),\quad t\in(0,1),
		\]
		which solves 
		\begin{equation}\label{hpequation}
			\begin{cases}
				\hpm''=(L^2\lambda_1(\omega)-p|\upm|^{p-1})\hpm \quad \mbox{ in } (0,1)\\
				\hpm(1)=-(\upm)'(1)\\
				\hpm'(0)=0.
			\end{cases}
		\end{equation}
		Being 
		\begin{equation}\label{rescalings}
			\lambda_1\left(\frac{\omega}{L}\right)=L^2\lambda_1(\omega)\qquad\text{and}\qquad\alpha_{1}(p,m)=L^2\alpha_{1,L}(p,m),
		\end{equation}
		\begin{equation}\label{iffs}
			\lambda_1(\omega)+\alpha_{1,L}(p,m)>0\quad \Leftrightarrow\quad \lambda_1\left(\frac{\omega}{L}\right)+\alpha_1(p,m)>0
			\quad \Leftrightarrow\quad L^2 \lambda_1(\omega)+\alpha_1(p,m)>0.
		\end{equation}
		Hence, if we assume $\lambda_1(\omega)+\alpha_{1,L}(p,m)>0$, then $\hpm$ is the unique solution to \eqref{hpequation} and $\hpLm'(L)\gtrless0$ if and only if $\hpm'(1)\gtrless 0$. Therefore the sign of $\hpm'(1)$ is what determines the stability/instability, by \emph{i)}, \emph{ii)} of Theorem \ref{teorema}.
		Furthermore by Theorem \ref{teorema}, \eqref{hpequation}, \eqref{rescalings} and \eqref{iffs} we have that $(\Omega_{\varphi_L}, \upLm)$ is a stable/unstable energy stationary pair if and only if $(\frac{\Omega_{\varphi_L}}{L},\upm)$ is a stable/unstable energy stationary pair.
	\end{osservazione}

	\subsection{Proof of Theorem \ref{Teo_stability_p_infty}}\label{subsection:proof_infty}
	In what follows, in order to simplify the notation, we set
	\[
	\lambda:=L^2\lambda_1(\omega),
	\]
	thus $\hpm$ solves \eqref{hptilde} with $\lambda=L^2\lambda_1(\omega)$.
	Furthermore we still denote by $\hpm$ the even extension to $I$ of $\hpm$.

	Next we prove a simple lemma.
	\begin{lemma}\label{lemma:ptoinfty_new}
		Let $m\geq1$ and $\gamma>m-\frac12$. There exists $p_1=p_1(m,\gamma)>1$
		such that, for every $p\geq p_1$ and every $\lambda>\gamma^2p^2$, we have
		$\lambda>-\alpha_1(p,m)$. In particular, the solution $\hpm$ of
		\eqref{hptilde} is positive.
	\end{lemma}
	\begin{proof}
		By Remark \ref{rem:tildealpha} and (5) of Theorem
		\ref{Teo_p_to_inf_rescaled}, applied with $n=2m-1$,
		\[
		\lambda+\alpha_1(p,m)=\lambda+\tilde\alpha_1(p,2m-1)
		>p^2\left(\gamma^2-\left(m-\frac12\right)^2+o_p(1)\right),
		\]
		where $o_p(1)\to0$ as $p\to+\infty$ and depends only on $m$. Since
		$\gamma>m-\frac12$, the right-hand side is positive for $p\geq p_1(m,\gamma)$.
		The positivity of $\hpm$ then follows from Proposition \ref{prop:hpL>0}.
	\end{proof}
	The key result, to prove Theorem \ref{Teo_stability_p_infty}, is the following.

	\begin{proposizione}\label{prop:derivative_hp}
		Let $m\geq1$ and $\gamma>m-\frac12$. There exists
		$p_2=p_2(m,\gamma)\geq p_1(m,\gamma)$ such that, for every $p\geq p_2$ and
		every $\lambda\in[\gamma^2p^2,\,p\|\upm\|_\infty^{p-1})$, the positive
		solution $\hpm$ of \eqref{hptilde} satisfies $\hpm'(1)>0$.
	\end{proposizione}
	\begin{proof}
		Let us suppose by contradiction that there exist a sequence $p_j\to+\infty$ and a sequence $\lambda_j\in[\gamma^2p_j^2,p_j\|u^{(m)}_{p_j}\|_\infty^{p_j-1})$ such that $\hpjm'(1)\leq0$.\\
		By definition of $\mu_{p_j}$, see \eqref{mu_p}, by \eqref{migliore_esp_p_to_infty} and since, by Theorem \ref{teo:globale_any_n} (1) and Proposition \ref{prop:monotonia_L_infty} (i), $\|\vpjdmu\|_\infty^{p_j-1}=(2m-1)^2\|v_{p_j}\|_\infty^{p_j-1}$, we have
		\[
		\frac{2\gamma^2}{(2m-1)^2}\longleftarrow\frac{\gamma^2 p_j^2}{(2m-1)^2 p_j \|v_{p_j}\|_\infty^{p_j-1}}\leq\lambda_j\frac{\mu_{p_j}^2}{(2m-1)^2}<\frac{p_j\|\vpjdmu\|_{\infty}^{p_j-1}}{ (2m-1)^2p_j\|v_{p_j}\|_\infty^{p_j-1}} = 1,
		\]
		and thus, up to a subsequence, \begin{equation}\label{NAnew}\mu_{p_j}^2\frac{\lambda_j}{(2m-1)^2}=\frac{{\lambda_j}}{(2m-1)^2p_j\|v_{p_j}\|_\infty^{p_j-1}}\to \eta \in \left[ \frac{2\gamma^2}{(2m-1)^2}, 1\right] \subset \left(\frac12,1 \right]\quad\text{ as $j \to +\infty$.}
		\end{equation}
		We divide the rest of the proof in two steps.\\
		For simplicity, in what follows we set
		\[
		r_i:=r_{2m-1,i},\qquad s_i:=s_{2m-1,i},
		\]
		for $i=1,\ldots,2m-1$, with $r_0=-1$.
		
		\emph{Step 1.} We show that the absolute maximum points of $h_{p_j,m}$ in $[r_{i-1},r_i]$ can not lie in $(r_{i-1},r_i)$.\\
		Let us fix $i\in\{1,\ldots,2m-1\}$ and let us suppose by contradiction that there exists $\xi_{p_j,m,i}\in(r_{i-1},r_i)$ such that
		\[
		h_{p_j,m}(\xi_{p_j,m,i}) =
		\max_{[r_{i-1},r_i]}h_{p_j,m}.
		\]
		Let us consider the rescaled function
		\[
		\kappa_{p_j,m,i}(y)=
		\frac{h_{p_j,m}\left(y\frac{\mu_{p_j}}{2m-1}+s_i\right)}
		{h_{p_j,m}(\xi_{p_j,m,i})},
		\quad
		y\in\left(-\frac{1}{\mu_{p_j}},\frac{1}{\mu_{p_j}}\right).
		\]
		If we set
		\[
		\tilde \xi_{p_j,m,i}
		:=(\xi_{p_j,m,i}-s_i)\frac{2m-1}{\mu_{p_j}}
		\in \left(-\frac{1}{\mu_{p_j}},\frac{1}{\mu_{p_j}}\right),
		\]
		then
		\begin{equation}\label{kappa_maximum}
			\kappa_{p_j,m,i}(\tilde\xi_{p_j,m,i})
			=
			1
			=
			\|\kappa_{p_j,m,i}\|_\infty
		\end{equation}
		and $\kappa_{p_j,m,i}$ solves the problem
		\begin{equation}\label{NAeq}
			\begin{cases}
				\begin{array}{ll}
					\tilde{\kappa}_{p_j,m,i}''(y)=
					\left(\lambda_j\frac{\mu_{p_j}^2}{(2m-1)^2}-       \left(1+\frac{\tilde{v}_{p_j}(y)}{p_j}\right)^{p_j-1}
					\right)\tilde{\kappa}_{p_j,m,i}(y)
					& y\in \left(-\frac{1}{\mu_{p_j}},\frac{1}{\mu_{p_j}}\right),
					\\
					\tilde{\kappa}_{p_j,m,i}>0
					&y\in \left(-\frac{1}{\mu_{p_j}},\frac{1}{\mu_{p_j}}\right),
					\\
					\|\tilde{\kappa}_{p_j,m,i}\|_{\infty}=1.&
				\end{array}
			\end{cases}
		\end{equation}
		where $\tilde v_{p_j}$ is defined in \eqref{utildep}.\\
		Being $\tilde v_{p_j}\leq0$ and by virtue of \eqref{NAnew},
		\begin{eqnarray*}
			|\kappa'_{p_j,m,i}(y)|
			&\leq&
			\left|
			\int_0^y\left|
			\lambda_j\frac{\mu_{p_j}^2}{(2m-1)^2}-
			\left(
			1+\frac{\tilde{v}_{p_j}(t)}{p_j}
			\right)^{p_j-1}
			\right|
			\kappa_{p_j,m,i}(t)\,dt
			\right|
			\\
			&\leq&
			\left|\int_0^y(1+1)\,dt\right|
			\leq2|y|,
		\end{eqnarray*}
		moreover by \eqref{NAeq}
		\[
		|\kappa''_{p_j,m,i}(y)|
		\leq
		\lambda_j\frac{\mu_{p_j}^2}{(2m-1)^2}+1
		\leq2.
		\]
		Thus $\kappa_{p_j,m,i}\to K_{m,i}$ in $C^1_{loc}(\R)$, where, by \eqref{NAnew}, $K_{m,i}$ is a solution to
		\begin{equation*}
			\begin{cases}
				\begin{array}{ll}
					K''=(\eta-e^W)K & \text{in $\R$,}\\
					K\geq0,\\
					\|K\|_\infty\leq1.
				\end{array}
			\end{cases}
		\end{equation*}
		Let us next show that $\|K_{m,i}\|_{\infty}=1$. Being $\tilde\xi_{p_j,m,i}$ a maximum point for $\kappa_{p_j,m,i}$, we have
		\begin{equation}\label{kappa_derivatives}
			\kappa'_{p_j,m,i}(\tilde\xi_{p_j,m,i})=0
			\quad\text{and}\quad
			\kappa''_{p_j,m,i}(\tilde\xi_{p_j,m,i})\leq0.
		\end{equation}
		If $|\tilde\xi_{p_j,m,i}|\to+\infty$, then by \eqref{kappa_derivatives}, \eqref{NAeq}, \eqref{kappa_maximum}, \eqref{NAnew} and \eqref{item:global_estimate_vtildep} of Proposition \ref{Teo_p_to_inf},
		\begin{eqnarray*}
			0
			&\geq&
			\kappa''_{p_j,m,i}(\tilde  \xi_{p_j,m,i})
			\\
			&=&
			\left(
			\lambda_j\frac{\mu_{p_j}^2}{(2m-1)^2}
			-\left(
			1+\frac{\tilde{v}_{p_j}(\tilde \xi_{p_j,m,i})}{p_j}
			\right)^{p_j-1}
			\right)
			\tilde{\kappa}_{p_j,m,i}(\tilde \xi_{p_j,m,i})
			\\
			&\geq&
			\eta+o_j(1)
			-
			e^{\frac{p_j-1}{p_j}\left(C_1 W(\tilde\xi_{p_j,m,i})+C_2\right)}
			\\
			&=&
			\eta+o_j(1),
		\end{eqnarray*}
		that is impossible. So we deduce that
		$|\tilde\xi_{p_j,m,i}|\leq C$ and, up to subsequences,
		$\tilde \xi_{p_j,m,i}\to\tilde \xi_{m,i}\in\R$, and
		\[
		1=\kappa_{p_j,m,i}(\tilde\xi_{p_j,m,i})
		\to
		K_{m,i}(\tilde\xi_{m,i}),
		\]
		which implies that $\|K_{m,i}\|_\infty=1$. From this we can also deduce that $K_{m,i}>0$.
		Summing up, $K_{m,i}$ solves \eqref{eq:tildeK}, so applying Lemma \ref{lemma:tildeK} we have that $\eta=\frac12$, which is a contradiction against $\eta\in(\frac12,1]$. This concludes Step 1.\\
		
		\emph{Step 2.} We prove that
		\[
		h'_{p_j,m}(r_i)>0,
		\qquad
		\forall\,i\in\{m,\ldots,2m-1\}.
		\]
		We proceed by induction on $i$.\\
		If $i=m$, we have $0\in(r_{m-1},r_m)$. By assumption and since $|\upjm|$ is strictly decreasing on $[0,r_m]$, $\lambda_j-p_j|\upjm|^{p_j-1}$ changes sign once in $[0,r_m]$. Therefore, $\hpjm''$ changes sign exactly once in $[0,r_m]$.\\
		Being $\hpjm'(0)=0$, $\hpjm''(0)<0$ and, by the evenness of $\hpjm$ in $[r_{m-1},r_m]$, if $\hpjm'(r_m)\leq0$, then $0$ is an absolute maximum point of $\hpjm$ in $[r_{m-1},r_m]$, but this is impossible in view of Step 1. Thus
		\[
		\hpjm'(r_m)>0.
		\]
		Now let us suppose that, for $i\in\{m+1,\ldots,2m-1\}$,
		\[
		\hpjm'(r_{i-1})>0.
		\]
		Then $r_{i-1}$ cannot be a maximum point of $\hpjm$ in $[r_{i-1},r_i]$; thus, by Step 1, $r_i$ is a maximum point of $\hpjm$ in $[r_{i-1},r_i]$, so in particular $\hpjm'(r_i)\geq0$. Finally, since by the equation $\hpjm''>0$ in a left neighborhood of $r_i$, the case $\hpjm'(r_i)=0$ can be ruled out and we conclude that $\hpjm'(r_i)>0$.\\
		Clearly Step 2, with $i=2m-1$, gives a contradiction against $\hpjm'(1)\leq0$, and this concludes the proof.
	\end{proof}

	We are finally in position to prove Theorem \ref{Teo_stability_p_infty}.
	
	\begin{proof} [Proof of Theorem \ref{Teo_stability_p_infty}]
		Let us first consider $(\Omega_{\varphi_L},\upLm)$. Let $p_0:=p_2(m,\gamma)$, where $p_2(m,\gamma)$ is given by Proposition \ref{prop:derivative_hp}, let $p\geq p_0$ and let $L>\gamma p/\sqrt{\lambda_1(\omega)}$, so that $\lambda=L^2\lambda_1(\omega)>\gamma^2p^2$. By Lemma \ref{lemma:ptoinfty_new}, \eqref{iffs} and Remark \ref{rem:nondegeneracy}, $\upLm$ is nondegenerate and the solution $\hpm$ of \eqref{hptilde} is positive.\\
		By Theorem \ref{teorema} and Remark \ref{rem:equivalence} it is sufficient to prove that $\hpm'(1)>0$.\\
		If $\lambda \geq p\|\upm\|_\infty^{p-1}$ then, for any $p>1$, by \eqref{hptilde} we get immediately that $\hpm''>0$ in $(0,1)$, which implies $\hpm'(1)>\hpm'(0)=0$.\\
		If $\lambda\in [{\gamma^2 p^2},p\|\upm\|_\infty^{p-1})$, then we can apply Proposition $\ref{prop:derivative_hp}$, obtaining that $\hpm'(1)>0$ and this shows the stability of $(\Omega_{\varphi_L},\upLm)$.\\
		By Remark \ref{rem:casonegativo} we deduce the thesis also for $(\Omega_{\varphi_L},-\upLm)$.
	\end{proof}

	\subsection{Proof of Theorem \ref{Teo_stability_p_1}}\label{subsection:proof_1}
	As in the previous Subsection, we set
	\[
	\lambda:=L^2\lambda_1(\omega),
	\]
	thus $\hpm$ solves \eqref{hptilde} with $\lambda=L^2\lambda_1(\omega)$.\\
	We start proving a simple lemma.
	\begin{lemma}\label{lemma:pto1_new}
		Let $m\geq1$ and $\lambda>(m-1)m\pi^2$. There exists $p_0=p_0(m,\lambda)>1$ such that, for every $p\in(1,p_0)$, we have $\lambda>-\alpha_1(p,m)$. In particular, the solution $\hpm$ of $\eqref{hptilde}$ is positive.
	\end{lemma}
	\begin{proof}
		Being $L>\pi\sqrt{\frac{m^2-m}{\lambda_1(\omega)}}$, by \eqref{rescalings}, Remark \ref{rem:tildealpha} and $\eqref{converg_alpha1p}$ of Theorem \ref{Teo_p_to_1_rescaled}
		\[
		\lambda_1(\omega)+\alpha_{1,L}(p,m)=\lambda_1(\omega)+\frac{\alpha_{1}(p,m)}{L^2}\rightarrow \lambda_1(\omega)+\frac{(m-m^2)\pi^2}{L^2} >0\quad\text{as $p\to1$}.
		\] 
		Thus, for $p$ sufficiently close to $1$, \begin{equation}\label{fiore}\lambda_1(\omega)>-\alpha_{1,L}(p,m),
		\end{equation} 
		which, by \eqref{iffs}, is equivalent to $\lambda>-\alpha_1(p,m)$.
		By \eqref{fiore}, by virtue of \eqref{iffs} and Proposition \ref{prop:hpL>0} we have that the solution $\hpm$ of \eqref{hptilde} is positive for $p$ sufficiently close to 1.
	\end{proof}
	The key result, to prove Theorem \ref{Teo_stability_p_1}, is the following.
	\begin{proposizione}\label{prop:hp_1}
		Let $m\geq1$ and let $\lambda \in \left( (m-1)m\pi^2, (m-\frac{1}{2})^2\pi^2\right)$. Then there exists $p_0=p_0(m,\lambda)>1$ such that, for every $p\in(1,p_0)$, $\hpm'(1)<0$, where $\hpm$ is the positive solution of $\eqref{hptilde}$.
	\end{proposizione}
	
	\begin{proof}
		Recalling that $v_p^{(2m-1)}$ is the even extension of $\upm$ to $I$, by (3) of Theorem \ref{Teo_p_to_1_rescaled}, applied with $n=2m-1$, we know that
		\[
		p|v_p^{(2m-1)}(t)|^{p-1}
		\longrightarrow
		\left(m-\frac12\right)^2\pi^2
		\]
		locally uniformly in $I\setminus
		\{r_{2m-1,1},\ldots,r_{2m-1,2m-2}\}$.
		In particular,
		\[
		p|v_p^{(2m-1)}(t)|^{p-1}
		\longrightarrow
		\left(m-\frac12\right)^2\pi^2
		\qquad\text{for a.e. }t\in I.
		\]
		
		Moreover, by (1) of Theorem \ref{Teo_p_to_1_rescaled}, applied with $n=2m-1$,
		\[
		\left\|p|v_p^{(2m-1)}|^{p-1}\right\|_\infty
		\longrightarrow
		\left(m-\frac12\right)^2\pi^2,
		\qquad\text{as $p\to1$,}
		\]
		thus it follows that, for $p$ sufficiently close to $1$,
		\begin{equation}\label{stima_unif_ap}
			0\leq p|v_p^{(2m-1)}(t)|^{p-1}
			\leq 2\left(m-\frac12\right)^2\pi^2
			\qquad\text{for every }t\in I.
		\end{equation}
		Hence, by the dominated convergence theorem,
		\begin{equation}\label{ap_to_a}
			p|v_p^{(2m-1)}|^{p-1}
			\longrightarrow\left(m-\frac12\right)^2\pi^2
			\qquad\text{in }L^1(I).
		\end{equation}
		
		We normalize $\hpm$, by setting
		\[
		\bar h_{p,m}(t):=
		\frac{\hpm(t)}{\hpm(0)},
		\qquad t\in[0,1],
		\]
		which solves
		\[
		\begin{cases}
			\bar h_{p,m}''=\left(\lambda-p|v_p^{(2m-1)}|^{p-1}\right)\bar h_{p,m}
			&\text{in }(0,1),\\
			\bar h_{p,m}(0)=1,\\
			\bar h_{p,m}'(0)=0.
		\end{cases}
		\]
		
		Since $\lambda<\left(m-\frac12\right)^2\pi^2$, we can define
		\[
		H_m(t):=\cos\left(\sqrt{\left(m-\frac12\right)^2\pi^2-\lambda}\:t
		\right),
		\]
		which is the unique solution of the Cauchy problem
		\[
		\begin{cases}
			H_m''=\left(\lambda-\left(m-\frac12\right)^2\pi^2\right)H_m
			&\text{in }(0,1),\\
			H_m(0)=1,\\
			H_m'(0)=0.
		\end{cases}
		\]
		We claim that
		\begin{equation}\label{hpm_to_Hm}
			\bar h_{p,m}\longrightarrow H_m
			\qquad\text{in }C^1([0,1]),\qquad\text{as $p\to1^+$.}
		\end{equation}
		
		To prove this, set $
		z_{p,m}:=\bar h_{p,m}-H_m$. Then $z_{p,m}(0)=z_{p,m}'(0)=0$ 
		and, subtracting the equations satisfied by $\bar h_{p,m}$ and $H_m$, we obtain
		\[
		z_{p,m}''=
		\left(\lambda-p|v_p^{(2m-1)}|^{p-1}\right)z_{p,m}+\left(
		\left(m-\frac12\right)^2\pi^2
		-p|v_p^{(2m-1)}|^{p-1}
		\right)H_m.
		\]
		
		By \eqref{stima_unif_ap} there exists a constant $M_m>0$ such that, for every $p$ sufficiently close to $1$,
		\[
		\left|
		\lambda-p|v_p^{(2m-1)}(t)|^{p-1}
		\right|
		\leq M_m
		\qquad\text{for every $t\in[0,1]$.}
		\]
		Furthermore, for every $t\in[0,1]$, using $z_{p,m}'(0)=0$ and $\|H_m\|_\infty=1$, we get
		\[
		|z_{p,m}'(t)|=\left|
		\int_0^t z_{p,m}''(s)\,ds
		\right|\leq M_m\int_0^t|z_{p,m}(s)|\,ds+\int_0^t
		\left|p|v_p^{(2m-1)}(s)|^{p-1}-
		\left(m-\frac12\right)^2\pi^2
		\right|\,ds.
		\]
		On the other hand, since $z_{p,m}(0)=0$,
		\[
		|z_{p,m}(t)|=
		\left|\int_0^t z_{p,m}'(s)\,ds
		\right|\leq \int_0^t|z_{p,m}'(s)|\,ds.
		\]
		Therefore, setting $
		Y_{p,m}(t):=|z_{p,m}(t)|+|z_{p,m}'(t)|$,
		we obtain
		\[
		Y_{p,m}(t)
		\leq
		(M_m+1)\int_0^tY_{p,m}(s)\,ds+
		\int_0^1\left|
		p|v_p^{(2m-1)}(s)|^{p-1}
		-\left(m-\frac12\right)^2\pi^2
		\right|\,ds
		\qquad\text{for every }t\in[0,1].
		\]
		By Gronwall's inequality,
		\[
		Y_{p,m}(t)\leq e^{M_m+1}
		\int_0^1\left|
		p|v_p^{(2m-1)}(s)|^{p-1}
		-\left(m-\frac12\right)^2\pi^2
		\right|\,ds
		\qquad\text{for every }t\in[0,1].
		\]
		Consequently, by \eqref{ap_to_a},
		$
		\|z_{p,m}\|_{C^1([0,1])}\leq
		2\sup_{t\in[0,1]}Y_{p,m}(t)
		=o_p(1)$,
		which proves \eqref{hpm_to_Hm}.

		At last, since
		\[
		\lambda\in
		\left(
		m(m-1)\pi^2,\left(m-\frac12\right)^2\pi^2
		\right),
		\]
		we have
		\[
		0<\left(m-\frac12\right)^2\pi^2-\lambda
		<\left(m-\frac12\right)^2\pi^2-m(m-1)\pi^2=\frac{\pi^2}{4},
		\]
		hence
		\begin{equation}\label{H_m'(1)}
			H_m'(1)=-\sqrt{
				\left(m-\frac12\right)^2\pi^2-\lambda
			}\:\sin\left(
			\sqrt{\left(m-\frac12\right)^2\pi^2-\lambda}\right)<0.
		\end{equation}
		
		Finally, combining \eqref{hpm_to_Hm} and \eqref{H_m'(1)}, we obtain, for $p$ sufficiently close to $1$,
		\[
		\hpm'(1)=\hpm(0)\bar h_{p,m}'(1)<0.
		\]
		This concludes the proof.
	\end{proof}
	
	\

We are finally in position to prove Theorem \ref{Teo_stability_p_1}.
\begin{proof}[Proof of Theorem $\ref{Teo_stability_p_1}$] 
	By Lemma \ref{lemma:pto1_new}, \eqref{iffs} and Remark \ref{rem:nondegeneracy}, there exists $p_0=p_0(m,\lambda)>1$ such that, for every $p\in(1,p_0)$, $\upLm$ is nondegenerate and the solution $\hpm$ of \eqref{hptilde} is positive.\\
	Furthermore we can apply Theorem \ref{teorema} and so by Remark \ref{rem:equivalence}, in order to prove the stability/instability of the couple $({\Omega_{\varphi_L}}, \upLm)$, it suffices to show that:
	\begin{itemize}
		\item[${(i)}$] if $\lambda>(m-\frac12)^2\pi^2$ and $p$ is sufficiently close to $1$, then $\hpm'(1)>0$,
		\item[${(ii)}$] if $\lambda\in((m^2-m)\pi^2,(m-\frac12)^2\pi^2)$ and $p$ is sufficiently close to $1$, then $\hpm'(1)<0$.
	\end{itemize}
	Proposition \ref{prop:hp_1} immediately yields $(ii)$, whereas $(i)$ can be easily shown. Indeed if $\lambda > (m-\frac12)^2\pi^2$, then, recalling that $v_p^{(2m-1)}$ is the even extension of $\upm$ to $I$, by (1) of Theorem \ref{Teo_p_to_1_rescaled}, applied with $n=2m-1$, we have \[
	p\|\upm\|_\infty^{p-1}=
	p\|v_p^{(2m-1)}\|_\infty^{p-1}
	\longrightarrow
	\left(m-\frac12\right)^2\pi^2.
	\]
	Hence, for $p$ sufficiently close to $1$,
	\[
	\lambda>p\|\upm\|_\infty^{p-1}\geq p|\upm|^{p-1}
	\qquad\text{in }[0,1].
	\]
	Thus, by \eqref{hptilde}, $\hpm''(t)>0$ for all $t\in[0,1]$, which implies
	$\hpm'(1)>\hpm'(0)=0$.
	
	By Remark \ref{rem:casonegativo} we deduce the thesis also for $(\Omega_{\varphi_L},-\upLm)$.
\end{proof}

\end{document}